\documentclass[12pt,leqno]{article}
\usepackage{amsfonts}
\usepackage{amsmath, amsthm, amsfonts, amssymb, color}
\usepackage{mathrsfs}
\usepackage{url}
\usepackage{color}
\theoremstyle{definition}

\newcommand{\scr}[1]{\mathscr #1}
\definecolor{wco}{rgb}{0.5,0.2,0.3}

\numberwithin{equation}{section} \theoremstyle{remark}

\newcommand{\ua}{\uparrow}

\title{{\bf Convergence  in  Wasserstein Distance for Empirical Measures of Dirichlet Diffusion Processes on Manifolds}\footnote{Supported in
 part by  NNSFC (11771326, 11831014, 11921001).} }
\author{
{\bf    Feng-Yu Wang$^{a),b)}$    }\\
\footnotesize{$^{a)}$ Center for Applied Mathematics, Tianjin University, Tianjin 300072, China }\\
 \footnotesize{ $^{b)}$ Department of Mathematics,
Swansea University,
Bay Campus,
Swansea,
SA1 8EN, United Kingdom}  }

\begin{document}
\allowdisplaybreaks
\def\R{\mathbb R}  \def\ff{\frac} \def\ss{\sqrt} \def\B{\mathbf
B}\def\TO{\mathbb T}
\def\I{\mathbb I_{\pp M}}\def\p<{\preceq}
\def\N{\mathbb N} \def\kk{\kappa} \def\m{{\bf m}}
\def\ee{\varepsilon}\def\ddd{D^*}
\def\dd{\delta} \def\DD{\Delta} \def\vv{\varepsilon} \def\rr{\rho}
\def\<{\langle} \def\>{\rangle} \def\GG{\Gamma} \def\gg{\gamma}
  \def\nn{\nabla} \def\pp{\partial} \def\E{\mathbb E}
\def\d{\text{\rm{d}}} \def\bb{\beta} \def\aa{\alpha} \def\D{\scr D}
  \def\si{\sigma} \def\ess{\text{\rm{ess}}}
\def\beg{\begin} \def\beq{\begin{equation}}  \def\F{\scr F}
\def\Ric{{\rm Ric}} \def\Hess{\text{\rm{Hess}}}
\def\e{\text{\rm{e}}} \def\ua{\underline a} \def\OO{\Omega}  \def\oo{\omega}
 \def\tt{\tilde}
\def\cut{\text{\rm{cut}}} \def\P{\mathbb P} \def\ifn{I_n(f^{\bigotimes n})}
\def\C{\scr C}      \def\aaa{\mathbf{r}}     \def\r{r}
\def\gap{\text{\rm{gap}}} \def\prr{\pi_{{\bf m},\varrho}}  \def\r{\mathbf r}
\def\Z{\mathbb Z} \def\vrr{\varrho} \def\ll{\lambda}
\def\L{\scr L}\def\Tt{\tt} \def\TT{\tt}\def\II{\mathbb I}
\def\i{{\rm in}}\def\Sect{{\rm Sect}}  \def\H{\mathbb H}
\def\M{\scr M}\def\Q{\mathbb Q} \def\texto{\text{o}} \def\LL{\Lambda}
\def\Rank{{\rm Rank}} \def\B{\scr B} \def\i{{\rm i}} \def\HR{\hat{\R}^d}
\def\to{\rightarrow}\def\l{\ell}\def\iint{\int}
\def\EE{\scr E}\def\Cut{{\rm Cut}}\def\W{\mathbb W}
\def\A{\scr A} \def\Lip{{\rm Lip}}\def\S{\mathbb S}
\def\BB{\scr B}\def\Ent{{\rm Ent}} \def\i{{\rm i}}\def\itparallel{{\it\parallel}}
\def\g{{\mathbf g}}\def\Sect{{\mathcal Sec}}\def\T{\mathcal T}\def\V{{\bf V}}
\def\PP{{\bf P}}\def\HL{{\bf L}}\def\Id{{\rm Id}}\def\f{{\bf f}}\def\cut{{\rm cut}}

\def\BL{\scr A}

\maketitle

\begin{abstract} Let $M$ be a   $d$-dimensional connected compact Riemannian manifold   with    boundary $\partial M$, let $V\in C^2(M)$ such that $\mu({\rm d} x):={\rm e}^{V(x)}{\rm d} x$ is a probability measure, and let
$X_t$ be the diffusion process generated by $L:=\Delta+\nabla V$ with  $\tau:=\inf\{t\ge 0: X_t\in\partial M\}$. Consider the  empirical measure
 $\mu_t:=\frac 1 t \int_0^t \delta_{X_s}{\rm d} s$  under the condition $t<\tau$ for the diffusion process. If $d\le 3$, then for any initial distribution not fully supported on $\partial M$, 
 \begin{align*} &c\sum_{m=1}^\infty \frac{2}{(\lambda_m-\lambda_0)^2} \le  \liminf_{t\to \infty} \inf_{T\ge t} \Big\{t {\mathbb E}\big[\mathbb W_2(\mu_t, \mu_0)^2\big|T<\tau\big]\Big\} \\
&\le  \limsup_{t\to \infty} \sup_{T\ge t} \Big\{ t \mathbb E\big[\mathbb W_2(\mu_t, \mu_0)^2\big|T<\tau\big] \Big\}\le \sum_{m=1}^\infty \frac{2}{(\lambda_m-\lambda_0)^2}\end{align*}
holds for some constant $c\in (0,1]$  with  $c=1$ when $\partial M$ is convex, 
 where $\mu_0:= \phi_0^2\mu$ for the first Dirichet eigenfunction $\phi_0$ of $L$,   $\{\lambda_m\}_{m\ge 0}$ are the Dirichlet eigenvalues of $-L$ listed in the increasing order counting multiplicities, and the  upper bound is finite if and only if $d\le 3$. 
 When $d=4$, $\sup_{T\ge t} \mathbb E\big[\mathbb W_2(\mu_t, \mu_0)^2\big|T<\tau\big] $  decays   in the order $t^{-1}\log t$, while for $d\ge 5$ it behaves like  $t^{-\frac 2 {d-2}}$, as $t\to\infty$. 
  \end{abstract} \noindent
 \vskip 2cm

\section{Introduction }

Let $M$ be a $d$-dimensional    connected  complete Riemannian manifold    with a smooth  boundary $\pp M$.  Let $V\in C^2(M)$ such that $\mu(\d x)=\e^{V(x)}\d x$ is a probability measure on $M$,
where $\d x$ is the Riemannian volume measure.
Let $X_t$  be the diffusion process generated by $L:=\DD+\nn V$ with hitting time
$$\tau:=\inf\{t\ge 0: X_t\in\pp M\}.$$ Denote by $\scr P$   the set of all probability measures on $M$, and let $\E^\nu$ be the expectation taken for the diffusion process with initial distribution $\nu\in\scr P$.
We  consider the  empirical measure
$$\mu_t:= \ff 1 t \int_0^t \dd_{X_s}\d s,\ \ t>0$$ under the condition that $t<\tau$. 
Since $\tau=0$ when $X_0\in \pp M$, to ensure $\P^\nu(\tau>t)>0,$ where $\P^\nu$ is the probability taken for the diffusion process with initial distribution $\nu$, we only consider
$$\nu\in \scr P_0:=\big\{\nu\in \scr P:\ \nu(M^\circ)>0\big\},\ \ M^\circ:= M\setminus \pp M.$$
Let $\mu_0=\phi_0^2\mu$, where $\phi_0$ is the first Dirichlet eigenfunction.
We investigate the convergence rate of $\E^\nu [\W_2(\mu_t, \mu_0)^2|t<\tau]$ as $t\to\infty$, where $\W_2$ is 
     the $L^2$-Wasserstein distance   induced by the Riemannian metric $\rr$. In general, for any $p\ge 1$, 
 $$\W_p(\mu_1,\mu_2):= \inf_{\pi\in \C(\mu_1,\mu_2)} \bigg(\int_{M\times M} \rr(x,y)^p \pi(\d x,\d y) \bigg)^{\ff 1 p},\ \ \mu_1,\mu_2\in \scr P,$$
 where $\C(\mu_1,\mu_2)$ is the set of all probability measures on $M\times M$ with marginal distributions $\mu_1$ and $\mu_2$, and  $\rr(x,y)$ is the Riemannian distance between $x$ and $y$, i.e.   the length of the shortest curve on $M$ linking $x$ and $y$.

 Recently, the convergence rate under $\W_2$ has been characterized  in \cite{WZ20} for the empirical measures of the $L$-diffusion processes without boundary (i.e. $\pp M=\emptyset$) or with a reflecting boundary.
Moreover,   the convergence of $\W_2(\mu_t^\nu,\mu_0)$ for the conditional empirical measure 
$$\mu_t^\nu:= \E^\nu(\mu_t|t<\tau),\ \ t>0$$   is investigated in \cite{W20}. Comparing with  $\E^\nu[\W_2(\mu_t, \mu_0)^2|t<\tau]$, in $\mu_t^\nu$  the conditional expectation inside the Wasserstein distance. 
According to \cite{W20}, $\W_2(\mu_t^\nu,\mu_0)^2$ behaves as $t^{-2}$, whereas  the following result says that  $\E[\W_2(\mu_t, \mu_0)^2|t<\tau]$ decays   at a slower rate, which coincides with the rate  of $\E[\W_2(\hat \mu_t, \mu)^2]$ given by  \cite[Theorems 1.1, 1.2]{WZ20}, where $ \hat \mu_t$ is the empirical measure of the reflecting diffusion process generated by $L$. 
 
  \beg{thm}\label{T1.1} Let   $\{\ll_m\}_{m\ge 0}$ be the Dirichlet  eigenvalues of $-L$ listed in the increasing order counting multiplicities.
Then  for any  $\nu\in \scr P_0$, the following assertions hold.
  \beg{enumerate} \item[$(1)$] In general, 
\beq\label{R0} \limsup_{t\to\infty}\Big\{t\sup_{T\ge t}   \E^\nu\big[\W_2(\mu_{t}, \mu_0)^2\big|T<\tau\big]\Big\}\le   \sum_{m=1}^\infty \ff{ 2}{(\ll_m-\ll_0)^2},\end{equation} 
and there exists a constant $c>0$ such that 
\beq\label{R00} \liminf_{t\to\infty}\Big\{ t\inf_{T\ge t}   \E^\nu\big[\W_2(\mu_{t}, \mu_0)^2\big|T<\tau\big]\Big\}\ge c   \sum_{m=1}^\infty \ff{ 2}{(\ll_m-\ll_0)^2}. \end{equation} 
If $\pp M$ is convex, then    $\eqref{R00}$ holds for $c=1$ so that 
$$ \lim_{t\to\infty}\Big\{ t \E^\nu\big[\W_2(\mu_{t}, \mu_0)^2\big|T<\tau\big]\Big\}=   \sum_{m=1}^\infty \ff{ 2}{(\ll_m-\ll_0)^2}\ \text{uniformly\ in\ }T\ge t.$$  
\item[$(2)$] When $d=4$, there exists a constant $c>0$ such that 
 \beq\label{R3} \sup_{T\ge t}   \E^\nu\big[\W_2(\mu_{t}, \mu_0)^2\big|T<\tau\big]\le  c t^{-1} \log t,\ \ t\ge 2. \end{equation} 
 \item[$(3)$] When $d\ge 5$, there exist a constant $c>1$ such that 
 \beg{align*}  c^{-1}  t^{-\ff 2{d-2}}  \le   \E^\nu\big[\W_1(\mu_{t}, \mu_0)^2\big|T<\tau\big]
  \le     \E^\nu\big[\W_2(\mu_{t}, \mu_0)^2\big|T<\tau\big]\le  c t^{-\ff 2{d-2}},\ \ T\ge t\ge 2. \end{align*}  
 \end{enumerate}
\end{thm}

Let $X_t^0$ be the diffusion process generated by $L_0:= L+2\nn\log\phi_0$ in $M^\circ$. It is well known that for any initial distribution supported on $M^\circ$, 
the law of $\{X^0_s: s\in [0,t]\}$ is the weak limit of the conditional distribution of $\{X_s: s\in [0,t]\}$ given $T<\tau$ as $T\to\infty$. Therefore, the following is a direct consequence of Theorem \ref{T1.1}.

\beg{cor} Let $\mu_t^0=\ff 1 t \int_0^t \dd_{X_s^0}\d s.$ Let $\nu\in \scr P_0$ with   $\nu(M^\circ)=1$. 
   \beg{enumerate} \item[$(1)$] In general, 
$$  \limsup_{t\to\infty} \Big\{ t    \E^\nu\big[\W_2(\mu^0_{t}, \mu_0)^2\big]\Big\}\le   \sum_{m=1}^\infty \ff{ 2}{(\ll_m-\ll_0)^2}, $$
and there exists a constant $c>0$ such that 
$$\liminf_{t\to\infty}\Big\{ t\inf_{T\ge t}    \big[\W_2(\mu^0_{t}, \mu_0)^2\big]\Big\}\ge c   \sum_{m=1}^\infty \ff{ 2}{(\ll_m-\ll_0)^2}.$$
If $\pp M$ is convex, then
 $$ \lim_{t\to\infty} \Big\{t   \E^\nu\big[\W_2(\mu_{t}, \mu_0)^2 \big]\Big\}=  \sum_{m=1}^\infty \ff{ 2}{(\ll_m-\ll_0)^2}.  $$
\item[$(2)$] When $d=4$, there exists a constant $c>0$ such that 
 $$  \E^\nu\big[\W_2(\mu^0_{t}, \mu_0)^2\big] \le  c t^{-1} \log t,\ \ t\ge 2.  $$
 \item[$(3)$] When $d\ge 5$, there exists a constant $c>1$ such that 
 $$ c^{-1}  t^{-\ff 2{d-2}} \le  \E^\nu\big[\W_2(\mu^0_{t}, \mu_0)^2\big] \le  c t^{-\ff 2{d-2}},\ \ t\ge 2.$$
  \end{enumerate}
\end{cor} 

In the next section, we first recall some facts on the Dirichlet semigroup and the diffusion semigroup $P_t^0$ generated by $L_0:= L+2\nn\log\phi_0$, then establish the Bismut derivative formula for $P_t^0$ which will be used to estimate the  lower bound of $\E^\nu[\W_2(\mu_t,\mu_0)^2|t<\tau)$. With these preparations, we prove Propositions \ref{Pn1} and \ref{PN2} in Sections 3 and 4 respectively, which imply Theorem \ref{T1.1}. 

 \section{Some preparations}

As in \cite{WZ20}, we first recall some well known facts on the Dirichlet semigroup, see for instances \cite{Chavel, Davies, OUB, W14}.
 Let $\{\phi_m\}_{m\ge 0}$ be the eigenbasis of the Dirichlet operator $L$ in $L^2(\mu)$, with Dirichlet eigenvalues $\{\ll_m\}_{m\ge 0}$ of $-L$ listed in the increasing order counting multiplicities. Then $\ll_0>0$ and
 \beq\label{EG}\|\phi_m\|_\infty\le \aa_0 \ss m,\ \  \aa_0^{-1} m^{\ff 2 d}\le \ll_m-\ll_0\le \aa_0 m^{\ff 2 d},\ \ m\ge 1\end{equation}
 holds for some constant $\aa_0>1$. Let $\rr_\pp$ be the Riemannian distance function to the boundary $\pp M$. Then $\phi_0^{-1}\rr_\pp$ is bounded such that
 \beq\label{JSS} \|\phi_0^{-1}\|_{L^p(\mu_0)}<\infty,\ \ p\in [1,3).\end{equation}
 The Dirichlet heat kernel has the representation
\beq\label{D09} p_t^D(x,y)=\sum_{m=0}^\infty \e^{-\ll_m t} \phi_m(x)\phi_m(y),\ \ t>0, x,y\in M.\end{equation}  Let $\E^x$ denote the expectation for the $L$-diffusion process starting at point $x$. Then Dirichlet diffusion semigroup generated by $L$ is given by
\beq\label{2.1} \beg{split} &P_t^D f(x):= \E^x[f(X_t)1_{\{t<\tau\}}]=\int_Mp_t^D(x,y)f(y)\mu(\d y)\\
&=\sum_{m=0}^\infty \e^{-\ll_m t} \mu(\phi_mf)\phi_m(x),\ \ t>0, f\in L^2(\mu).\end{split}\end{equation}
Consequently,  
\beq\label{D10} \lim_{t\to\infty} \big\{\e^{\ll_0 t} \P^\nu(t<\tau) \big\}=  \lim_{t\to\infty} \big\{\e^{\ll_0 t} \nu(P_t^D1) \big\} =\mu(\phi_0)\nu(\phi_0),\ \ \nu\in \scr P_0.\end{equation} 
Moreover, there exists a constant $c>0$ such that
 \beq\label{AC0} \|P_t^D\|_{L^p(\mu)\to L^q(\mu)}:= \sup_{\mu(|f|^p)\le 1} \|P_t^Df\|_{L^q(\mu)}\le c\e^{-\ll_0t} (1\land t)^{-\ff {d (q-p)}{2pq}}, \ \ t>0, q\ge p\ge 1.\end{equation}

On the other hand,  let $L_0= L+2\nn\log \phi_0$. Noting that $L_0f= \phi_0^{-1} L(f\phi_0) + \ll_0 f$,  $L_0$ is a self-adjoint operator in $L^2(\mu_0)$ and the associated semigroup $P_t^0:=\e^{t L_0}$ satisfies 
 \beq\label{PR0} P_t^0f=\e^{\ll_0 t}\phi_0^{-1} P_t^D(f\phi_0),\ \ f\in L^2(\mu_0),\ \ t\ge 0.\end{equation}
So,    $\{\phi_0^{-1}\phi_m\}_{m\ge 0}$ is an eigenbasis of $L_0$ in $L^2(\mu_0)$ with
 \beq\label{PR1} L_0(\phi_m\phi_0^{-1})= -(\ll_m-\ll_0)\phi_m\phi_0^{-1},\ \ P_t^0(\phi_m\phi_0^{-1})=\e^{-(\ll_m-\ll_0)t} \phi_m\phi_0^{-1},\ \ m\ge 0, t\ge 0.\end{equation}
 Consequently,
 \beq\label{ONN} P_t^0f = \sum_{m=0}^\infty \mu_0(f\phi_m\phi_0^{-1}) \e^{-(\ll_m-\ll_0)t} \phi_m\phi_0^{-1},\ \ f\in L^2(\mu_0),\end{equation}
 and the heat kernel of $P_t^0$ with respect to $\mu_0$ is given by
 \beq\label{ON*} p_t^0(x,y)= \sum_{m=0}^\infty (\phi_m\phi_0^{-1})(x)  (\phi_m\phi_0^{-1})(y) \e^{-(\ll_m-\ll_0)t},\ \ x,y\in M, t>0.\end{equation}
 By the intrinsic ultracontractivity, see for instance \cite{07OW}, there exists a constant $\aa_1\ge 1$ such that
\beq\label{PRR0}  \|P_t^0-\mu_0\|_{L^1(\mu_0)\to L^\infty(\mu_0)}:= \sup_{\mu_0(|f|)\le 1} \|P_t^0f-\mu_0(f)\|_{\infty} \le \ff {\aa_1\e^{-(\ll_1-\ll_0)t}}{(1\land t)^{\ff{d+2}2}},\ \ t>0.\end{equation}
Combining this with the semigroup property and the contraction of $P_t^0$ in $L^p(\mu)$ for any $p\ge 1$, we find a constant $\aa_2\ge 1$ such that
\beq\label{000} \|P_t^0-\mu_0\|_{L^p(\mu_0)}:=  \sup_{\mu_0(|f|^p)\le 1} \|P_t^0f-\mu_0(f)\|_{L^p(\mu_0)} \le \aa_2\e^{-(\ll_1-\ll_0)t},\ \ t\ge 0, p\ge 1.\end{equation}
By the interpolation theorem,  \eqref{PRR0} and \eqref{000} yield that for some constant $\aa_3>0$, 
\beq\label{**1} \|P_t^0-\mu_0\|_{L^p(\mu_0)\to L^q(\mu_0)}\le  \aa_3 \e^{-(\ll_1-\ll_0)t} \{1\land t\}^{-\ff{(d+2)(q-p)}{2pq}},\ \ t>0,\infty\ge q>p\ge 1.\end{equation}
By this and \eqref{PR1},    there exists  a constant $\aa_4>0$ such that
\beq\label{CC} \|\phi_m\phi_0^{-1}\|_\infty \le \aa_4 m^{\ff{d+2}{2d}},\ \ m\ge 1.\end{equation}

\

In the remainder of this section, we establish the Bismut derivative formula for   $P_t^0$, which is not included by existing results due to the singularity of $\nn\log \phi_0$ in $L_0$.    
 Let $X_t^0$ be the diffusion process generated by $L_0$, which solves the following It\^o SDE on $M^\circ$, see \cite{EM}:
\beq\label{E01} \d^I X_t^0= \nn (V+2\log\phi_0)(X_t^0)\d t + \ss 2 U_t \d B_t,\end{equation}
where $B_t$ is the $d$-dimensional Brownian motion, and $U_t\in O_{X_t^0}(M)$ is the horizontal lift of $X_t^0$ to the frame bundle $O(M)$. Let $\Ric$ and $\Hess$  be the Ricci curvature and the Hessian tensor on $M$ respectively. Then the Bakry-Emery curvature of $L_0$ is given by 
   $$\Ric_{L_0}  := \Ric -\Hess_{V+2\log \phi_0}.$$ Let $\Ric_{L^0}^\#(U_t)\in\R^d\otimes\R^d$ be defined by
  $$\<\Ric_{L^0}^\#(U_t) a,b\>_{\R^d} = \Ric_{L_0}(U_t a, U_t b),\ \ a,b\in \R^d.$$ 
 We consider the following ODE on $\R^d\otimes\R^d$:
  \beq\label{ODE} \ff{\d}{\d t} Q_t= -\Ric_{L^0}^\#(U_t) Q_t,\ \ Q_0=I,\end{equation}
  where $I$ is the identity matrix. 

\beg{lem}\label{L*} For any $\vv>0$, there exist  constants  $\dd_1,\dd_2>0$ such that 
\beq\label{KA0} \E^x\big[\e^{\dd_1 \int_0^t \{\phi_0(X_s)\}^{-2}\d s}\big] \le \dd_2 \phi_0^{-\vv} (x) \e^{\dd_2 t},\ \ t\ge 0, x\in M^\circ.\end{equation} 
Consequently, 
\beg{enumerate} \item[$(1)$] For any $\vv>0$ and $p>1$, there exists a constant $\kk >0$ such that 
$$|\nn P_tf(x)|^2\le \kk \phi_0(x)^{-\vv} \e^{\kk t} \{P_t|\nn f|^{2p}(x)\}^{\ff 1 p} ,\ \ f\in C_b^1(M).$$
\item[$(2)$] For any $\vv>0$ and $p\ge 1$, there exists a constant $\kk>0$ such that for any stopping time $\tau'$, 
$$ \E^x [\|Q_{t\land\tau'}\|^p]\le \kk \phi_0(x)^{-\vv} \e^{\kk t},\ \ t\ge 0.$$\end{enumerate} 
\end{lem} 

\beg{proof} Since  $L\phi_0=-\ll_0\phi_0$, $\phi_0>0$ in $M^\circ$,  $\|\phi_0\|_\infty<\infty$ and $|\nn \phi_0|$ is strictly positive in a neighborhood of $\pp M$, we find a constant $c_1,c_2>0$ such that 
$$L_0\log \phi_0^{-1} = -\phi_0^{-1} L\phi_0 + \phi_0^{-2}|\nn\phi_0|^2 - 2\phi_0^{-2} |\nn \phi_0|^2\le c_1 -c_2 \phi_0^{-2}.$$
So, by \eqref{E01} and It\^o's formula, we obtain 
$$\d \log \phi_0^{-1} (X_t^0) \le \{c_1-c_2\phi_0^{-2}(X_t^0)\}\d t + \ss 2 \<\nn\log \phi_0^{-1}(X_t^0), U_t\d B_t\>.$$ This implies
\beq\label{*B*} \E^x\int_0^t [\phi_0^{-2}(X_s^0)]\d s \le ct+ c \log (1+\phi_0^{-1})(x),\ \ t\ge 0 \end{equation}
for some constant $c>0$, and 
  for any constant $\dd>0$, 
\beg{align*} &\E^x\big[\e^{\dd c_2\int_0^t \phi_0^{-2}(X_s^0)\}\d s}\big]\le \E^x\big[\e^{\dd \log\phi_0^{-1}(x) + \dd \log \phi_0(X_t^0) + c_1 \dd t -\dd \ss 2 \int_0^t  \<\nn\log \phi_0(X_s^0), U_s\d B_s\>}\big]\\
&\le \e^{c_1\dd t} \phi_0^{-\dd} (x) \|\phi_0\|_\infty^{\dd} \big(\E^x[\e^{4\dd^2\int_0^t |\nn \log \phi_0|^2(X_s^0)\d s}]\big)^{\ff 1 2}.\end{align*}
Let $c_3= 4 \|\nn\phi_0\|^2_\infty$, and take $\dd\in (0, c_2/c_3]$, we derive
$$\E^x\big[\e^{\dd c_2\int_0^t \phi_0^{-2}(X_s^0)\}\d s}\big]\le \e^{2c_1\dd t} \phi_0^{-2\dd} (x),\ \ \dd\in (0, c_2/c_3].$$
 This implies \eqref{KA0}. 
Below we prove assertions (1) and (2) respectively.

 Since $V\in C_b^2(M)$ and $\phi_0\in C^2_b(M)$ with $\phi_0>0$ in $M^\circ$, there exists a constant $\aa_1>0$ such that 
\beq\label{CV}  \Ric_{L_0}(U,U)\ge -\aa_1\phi_0^{-1}(x) |U|^2,\ \ x\in M^\circ, U\in T_x M.\end{equation} 
By \eqref{E01},  \eqref{CV},  and  the formulas of It\^o and Bochner,  for fixed $t>0$ this implies 
\beg{align*} &\d |\nn P_{t-s}^0 f|^2(X_s^0) \\
&= \big\{L_0  |\nn P_{t-s}^0 f|^2(X_s^0)-2\<\nn P_{t-s}^0 f, \nn L_0 P_{t-s}^0 f\> \big\}\d s  +\ss 2 \<\nn  |\nn P_{t-s}^0 f|^2(X_s^0), U_s\d B_s\>\\
&\ge 2\Ric_{L^0}(\nn P_{t-s}^0 f,\nn P_{t-s}^0 f)(X_s^0) \d s + \ss 2 \<\nn  |\nn P_{t-s}^0 f|^2(X_s^0), U_s\d B_s\>\\
&\ge -2 \aa_1 \{\phi_0^{-1}|\nn P_{t-s}^0 f)|^2\}(X_s^0) \d s+ \ss 2 \<\nn  |\nn P_{t-s}^0 f|^2(X_s), U_s\d B_s\>\d s.\end{align*}
Then  
\beg{align*} &|\nn P_tf(x)|^2= \E^x |\nn P_t f|^2(X_0^0) \le   \E^x\big[|\nn  f|^2(X_t^0)\e^{2\int_0^t 2\aa_1 \phi^{-1}(X_u^0)\d u}\big]\\
&\le \big\{  \E^x \e^{\ff{2\aa_1p}{p-1} \int_0^t  \phi^{-1}(X_u^0)\d u}\big] \big\}^{\ff{p-1}p} \{P_t|\nn f|^{2p}(x)\}^{\ff 1 p}.\end{align*}
Combining this   with \eqref{KA0}, we prove   (1).

Next, by   \eqref{ODE} and  \eqref{CV}, we obtain 
$$ \|Q_{t\land\tau'}\|\le \e^{\aa_1\int_0^t \phi^{-1}(X_s^0)\d s},\ \ t\ge 0.$$
This together with \eqref{KA0}  implies  (2). 

\end{proof} 

\beg{lem}\label{LBS}  For any $t>0$ and $\gg \in C^1([0,t])$ with $\gg(0)=0$ and $\gg (t)=1$, we have 
\beq\label{BSF} \nn P_t^0 f(x)=  \E^x \bigg[f(X_t^0) \int_0^t \gg'(s) Q_s^*\d B_s\bigg],\ \ x\in M^\circ, f\in \B_b(M^\circ).\end{equation}
Consequently, for any $\vv>0$ and $p>1$, here exists a constant $c>0$ such that 
\beq\label{GRD} |\nn P_t^0f| \le \ff {c\phi_0^{-\vv}} {\ss {1\land t}} (P_t^0|f|^p)^{\ff 1 p},\ \ t>0,  f\in \B_b(M^\circ).\end{equation} 
\end{lem} 

\beg{proof} Since \eqref{GRD} follows from \eqref{BSF}  with $\gg(s):=\ff{t-s}t$ and Lemma \ref{L*}(2), it suffices to prove the Bismut formula  \eqref{BSF}.  By an approximation argument, we only need to prove for $f\in C_b^1(M)$. 
The proof is standard by Elworthy-Li's  martingale argument  \cite{EL}, see also \cite{Thalmaier}. 
By $\|\nn f\|_\infty<\infty$ and  Lemma \ref{L*}(1) for $\vv=\ff 1 4$,  we find a constant $c_1>0$ such that 
\beq\label{LB1} |\nn P_s^0 f|(x) \le c_1 \phi_0^{-1/4}(x),\ \ s\in [0,t], x\in M^\circ. \end{equation}

Next, since $L\phi_0=-\ll_0 \phi_0$ implies $L_0 \phi_0^{-1} = \ll_0 \phi_0^{-1}$, by It\^o's formula we obtain   
\beq\label{LB2} \E^x [\phi_0^{-1}(X_{t\land \tau_n}^0)] \le  \phi_0^{-1}(x) \e^{\ll_0 t},\ \ t\ge 0, n\ge 1,\end{equation}
where  $\tau_n:= \inf\{t\ge 0:\phi_0(X_s^0)\le \ff 1 n\}\uparrow \infty$  as $n\uparrow\infty$ by noting that the process $X_t^0$ is non-explosive in $M^\circ$.

Moreover,  by  It\^o's formula, for any $a\in \R^d$, we have 
\beg{align*} & \d\<\nn P_{t-s}^0f(X_s^0), U_s Q_sa\> =\ss 2\, \Hess_{P_{t-s}f} (U_s\d B_s, U_sQ_s a)(X_s^0),\\
& \d P_{t-s}f(X_s^0) =\ss 2\, \<\nn P_{t-s}^0 f(X_s^0), U_s\d B_s\>,\ \ s\in [0,t].\end{align*} 
Due to  the integration by part formula, this   and $\gg(0)=0$  imply 
\beq\label{KL} \beg{split} &- \ff 1 {\ss 2} \E^x\bigg[f(X_{t\land \tau_n}^0) \int_0^{t\land\tau_n}\gg'(s)  \<Q_s a, \d B_s\>\bigg] \\
&= \E\bigg[\int_0^{t\land\tau_n} \<\nn P_{t-s}^0 f(X_s^0), U_s Q_s a\>\d (1-\gg)(s)\bigg]\\  
&= \E \big[(1-\gg)(t\land\tau_n) \<\nn P_{t-t\land\tau_n}^0 f(X_{t\land \tau_n}^0), Q_{t\land \tau_n} a\> \big] - \<\nn P_t f(x), U_0 a\> \\
&\qquad  -\E\bigg[\int_0^{t\land\tau_n} (1-\gg)(s) \d \<\nn P_{t-s}^0 f(X_s^0), U_s Q_s a\> \bigg]\\ 
&= \E \big[(1-\gg)(t\land\tau_n) \<\nn P_{t-t\land\tau_n}^0 f(X_{t\land \tau_n}^0), Q_{t\land \tau_n} a\> \big] - \<\nn P_t f(x), U_0 a\>,\ \ n\ge 1.\end{split}\end{equation} 
Since $\gg$ is bounded with $\gg(t)=1$ such that $(1-\gg)(t\land\tau_n)\to 0$ as $n\to\infty,$ and     \eqref{LB1}, \eqref{LB2}  and Lemma \ref{L*}(2)  imply
$$\sup_{n\ge 1} \E^x \big[\<\nn P_{t-t\land\tau_n}^0 f(X_{t\land \tau_n}^0), Q_{t\land \tau_n} a\>^2\big] \le c_1 \sup_{n\ge 1} \big(\E[\phi_0^{-1}(X_{t\land\tau_n}^0)] \big)^{\ff 1 2} \big(\E^x \|Q_{t\land \tau_n}\|^{4}\big)^{\ff 1 2} <\infty,$$
by the dominated convergence theorem, we may take  $n\to\infty$ in  \eqref{KL} to derive \eqref{BSF}. 
\end{proof} 
  
 \section{Upper bound estimates}
 In this section we prove the following result which includes upper bound estimates in Theorem \ref{T1.1}.
 
 \beg{prp}\label{Pn1} Let    $\nu\in \scr P_0$. 
  \beg{enumerate} \item[$(1)$]   $\eqref{R0}$ holds. 
\item[$(2)$] When $d=4$, there exists a constant $c>0$ such that $\eqref{R3}$ holds. 
 \item[$(3)$] When $d\ge 5$, there exists a constant $c>0$ such that 
 $$ \sup_{T\ge t}   \E^\nu\big[\W_2(\mu_{t}, \mu_0)^2\big|T<\tau\big]\le  c t^{-\ff 2{d-2}},\ \ t\ge 2. $$
 \end{enumerate}
    \end{prp} 
The main tool in the study of the upper bound estimate is the following inequality due to \cite{AMB}, see also \cite[Lemma 2.3]{WZ20}: for any
probability density $g\in L^2(\mu_0)$,
\beq\label{*UPA} \W_2(g\mu_0,\mu_0)^2\le \int_M \ff{|\nn L_0(g-1)|^2}{\scr M(g,1)}\d \mu_0,\end{equation} 
where $\scr M(a,b):= \ff{a-b}{\log a -\log b}1_{\{a\land b>0\}}.$ To apply this inequality, as in \cite{WZ20}, we first modify $\mu_t$ by $\mu_{t,r} :=\mu_tP_r^0$ for some $r>0$,
where for a probability measure $\nu$ on $M^\circ$, $\nu P_r^0$ is the law of  the $L_0$-diffusion process $X_r^0$ with initial distribution $\nu$.   Obviously,  by \eqref{ON*} we have 
 \beq\label{DS1}\beg{split} & \rr_{t,r}:= \ff{\d \mu_{t,r}}{\d\mu_0} =    \ff 1 t \int_0^t p_r^0(X_s,\cdot)\d s= 1+ \sum_{m=1}^\infty \e^{-(\ll_m-\ll_0)r} \psi_m(t)   \phi_m\phi_0^{-1},\\
 & \psi_m(t):= \ff 1 t \int_0^t \{\phi_m\phi_0^{-1}\}(X_s)\d s,\end{split}\end{equation}
which are well-defined on the event $\{t<\tau\}$.

 \beg{lem}\label{LL1}        If  $d\le 3$ and $\nu= h\mu$ with $h\phi_0^{-1}\in L^p(\mu_0)$ for some $p>\ff{d+2} 2,$  then there exists a constant $c>0$ such that 
\beg{align*} & \sup_{T\ge t}   \Big|t  \E^\nu\big [ \mu_0(|\nn L_0^{-1} (\rr_{t,r}-1)|^2)\big|T<\tau\big]- 2\sum_{m=1}^\infty \ff{\e^{-2(\ll_m-\ll_0) r } }{(\ll_m-\ll_0)^2}\Big|\\
&\le c t^{-1} \big( r^{-\ff{(d-2)^+}2} +1_{\{d=2\}} \log r^{-1}\big),\ \ r\in (0,1], t\ge 1.\end{align*} 
\end{lem} 

 \beg{proof} By the Markov property,    \eqref{PR0} and \eqref{2.1}, we have 
 \beq\label{SD2}\beg{split} & \E^x [f(X_s)1_{\{T<\tau\}}] =   \E^x\big[1_{\{s<\tau\}} f(X_s) \E^{X_s} 1_{\{T-s<\tau\}}\big] \\
 &=   P_s^D\{fP_{T-s}^D 1\}(x) = \e^{-\ll_0 T}  \big( \phi_0 P_s^0 \{f P_{T-s}^0\phi_0^{-1}\}\big)(x),\ \ s<T.\end{split} \end{equation}
 By the same reason,   and noting that  $\E^\nu = \int_M \E^x \nu(\d x)$, we derive 
 \beg{align*} &  \E^\nu [f(X_{s_1})f(X_{s_2}))1_{\{T<\tau\}}] =  \int_M \E^x\big[1_{\{s_1<\tau\}} f(X_{s_1}) \E^{X_{s_1}} \{f(X_{s_2-s_1}) 1_{\{T-s_1<\tau\}}\}\big]\nu(\d x)\\
 &=\e^{-\ll_0T}  \nu\big(\phi_0P_{s_1}^0[fP_{s_2-s_1}^0\{fP_{T-s_2}^0\phi_0^{-1}\}]\big),\ \ s_1<s_2<T.\end{align*}In particular, the formula with 
 $f=1$ yields 
 $$\P^\nu(T<\tau)=\e^{-\ll_0 T}  \nu(\phi_0 P_T^0 \phi_0^{-1}).$$
Combining  these   with \eqref{DS1}, \eqref{PR1}, $\E^\nu(\xi|T<\tau):= \ff{\E^\nu [\xi 1_{\{T<\tau\}}]}{\P^\nu(T<\tau)}$ for an integrable random variable $\xi$, and the symmetry of $P_t^0$ in $L^2(\mu_0)$,  for $\nu=h\mu$  we obtain
 \beq\label{SD3} \beg{split} &t   \E^\nu\big[ \mu_0(|\nn L_0^{-1} (\rr_{t,r}-1)|^2)\big|T<\tau\big]= \sum_{m=1}^\infty \ff {t \E^\nu[\psi_m(t)^2|T<\tau]}{\e^{2(\ll_m-\ll_0)r}(\ll_m-\ll_0)} \\
&= \sum_{m=1}^\infty \ff {2\int_0^t  \d s_1\int_{s_1}^t \E^\nu\big[1_{\{T<\tau\}} (\phi_m\phi_0^{-1})(X_{s_1})(\phi_m\phi_0^{-1})(X_{s_2}) \big] \d s_2} {t\e^{2(\ll_m-\ll_0)r}(\ll_m-\ll_0)\nu(\phi_0P_T^0\phi_0^{-1})} \\
 & = \sum_{m=1}^\infty  \ff {2\int_0^t\d s_1\int_{s_1}^t \nu\big(\phi_0^{-1} P_{s_1}^0 \{\phi_m\phi_0^{-1} P_{s_2-s_1}^0 [\phi_m\phi_0^{-1} P_{T-s_2}^0 \phi_0^{-1}]\}\big)\d s_2} {t\e^{2(\ll_m-\ll_0)r}(\ll_m-\ll_0)\nu(\phi_0P_T^0\phi_0^{-1})}\\
 &= \sum_{m=1}^\infty \ff {2\int_0^t\d s_1\int_{s_1}^t \mu_0\big(\{P_{s_1}^0 (h \phi_0^{-1})\} \phi_m\phi_0^{-1} P_{s_2-s_1}^0 [\phi_m\phi_0^{-1} P_{T-s_2}^0 \phi_0^{-1}]\big)\d s_2} {t\e^{2(\ll_m-\ll_0)r}(\ll_m-\ll_0)\mu_0(\phi_0^{-1} P_T^0(h\phi_0^{-1}) )}.\end{split}\end{equation}
 By \eqref{**1}, $\|\phi_0^{-1}\|_{L^2(\mu_0)}=1$ and $\|h\phi_0^{-1}\|_{L^1(\mu_0)}=\mu(h\phi_0)\le \|\phi_0\|_\infty<\infty$, we find a constant $c_1>0$ such that 
 \beq\label{B_1} \beg{split} &\big|\mu_0(\phi_0^{-1} P_T^0(h\phi_0^{-1}))-\mu(\phi_0)\nu(\phi_0)|\le \|\phi_0^{-1} (P_T^0-\mu_0) (h\phi_0^{-1})\|_{L^1(\mu_0)}\\
 &\le \|P_T^0-\mu_0\|_{L^1(\mu_0)\to L^2(\mu_0)}  \|h\phi_0^{-1}\|_{L^1(\mu_0)}\le c_1 \e^{-(\ll_1-\ll_0)T},\ \ T\ge 1.\end{split}\end{equation}
 On the other hand, write 
  \beq\label{SD7}\beg{split} & \mu_0\big(\{P_{s_1}^0 (h \phi_0^{-1})\}    \phi_m\phi_0^{-1} P_{s_2-s_1}^0 [\phi_m\phi_0^{-1} P_{T-s_2}^0 \phi_0^{-1}]\big)\\
  &=\nu(\phi_0)\mu(\phi_0)\e^{-(\ll_m-\ll_0)(s_2-s_1)} + J_1(s_1,s_2)+J_2(s_1,s_2)+J_3(s_1,s_2),\end{split}\end{equation}
  where, due to \eqref{PR1}, 
 \beg{align*} &J_1(s_1,s_2):= \mu_0\big(\{P_{s_1}^0 (h \phi_0^{-1})-\mu(h\phi_0)\}  \phi_m\phi_0^{-1} P_{s_2-s_1}^0 [\phi_m\phi_0^{-1} (P_{T-s_2}^0 \phi_0^{-1}-\mu(\phi_0))]\big),\\
 &J_2(s_1,s_2):= \mu(\phi_0)\e^{-(\ll_m-\ll_0)(s_2-s_1)} \mu_0\big(\{P_{s_1}^0 (h \phi_0^{-1})-\mu(h\phi_0)\}    \{\phi_m\phi_0^{-1} \}^2 \big),\\
 &J_3(s_1,s_2):= \mu(h\phi_0)\e^{-(\ll_m-\ll_0)(s_2-s_1)} \mu_0\big(    \{\phi_m\phi_0^{-1} \}^2  \{P_{T-s_2}^0 \phi_0^{-1}]-\mu(\phi_0)\}\big).\end{align*}  
By \eqref{SD3}, \eqref{B_1} and \eqref{SD7}, we find a constant $\kk>0$ such that 
\beq\label{007} \beg{split} &\sup_{T\ge t}   \Big|t  \E^\nu\big [ \mu_0(|\nn L_0^{-1} (\rr_{t,r}-1)|^2)\big|T<\tau\big]- 2\sum_{m=1}^\infty \ff{\e^{-2(\ll_m-\ll_0) r } }{(\ll_m-\ll_0)^2}\Big|\\
&\le \ff \kk t  \sum_{m=1}^\infty\bigg(\ff{\e^{-2(\ll_m-\ll_0)r}}{(\ll_m-\ll_0)^2}+ \ff{\e^{-2(\ll_m-\ll_0)r}}{\ll_m-\ll_0} \int_0^t\d s_1\int_{s_1}^t |J_1+J_2+J_3|(s_2,s_2)\d s_2\bigg),\ \ t\ge 1.\end{split}\end{equation} 
Since    $\|h\phi_0^{-1}\|_{L^p(\mu_0)}<\infty$,  $\|\phi_0^{-1}\|_{L^\theta(\mu_0)}<\infty$ for $\theta<3$ due to \eqref{JSS},  $\|\phi_m\phi_0^{-1}\|_{L^2(\mu_0)}=1$,  by \eqref{**1}, for any $\theta\in (\ff 5 2,3)$, we find constants $c_1,c_2>0$ such that   
 \beq\label{61} \beg{split}  &|J_1(s_1,s_2)|\le c_1 \|P_{s_1}^0-\mu_0\|_{L^p(\mu_0)\to L^\infty(\mu_0)}    \|P_{T-s_2}^0 -\mu_0\|_{L^\theta(\mu_0)  \to L^\infty(\mu_0) } \\ 
   &\le c_2   \e^{-(\ll_1-\ll_0)(s_1+T-s_2)} (1\land s_1)^{-\ff{d+2}{2p}} \{1\land (T-s_2)\}^{-\ff{d+2}{ 2\theta} },\end{split} \end{equation}
and    
 \beq\label{62} \beg{split}  &|(J_2 +J_3)(s_1,s_2)|\\
 &\le c_1\e^{-(\ll_m-\ll_0)(s_2-s_1)}\big( \|P_{s_1}^0-\mu_0\|_{p\to \infty}  +  \|P_{T-s_2}^0-\mu_0\|_{L^\theta (\mu_0) \to L^\infty(\mu_0)}   \big)\\
 &\le c_2 \e^{-(\ll_m-\ll_0)(s_2-s_1)} \big(\{1\land s_1\}^{-\ff{d+2}{2p}} \e^{-(\ll_1-\ll_0)s_1}    
 +     \{1\land (T-s_2)\}^{-\ff{d+2}{2 \theta} }\e^{-(\ll_1-\ll_0)(t-s_2)}\big). \end{split}\end{equation}
Since $q>\ff 5 2 $ and $p>\ff{d+2}2$ imply $\ff{d+2}{2q} \lor  \ff{d+2}{2p}<1$ for $d\le 3$,  by \eqref{61} and \eqref{62}, 
we find a constant $c>0$ such that  
$$\int_0^t\d s_1\int_{s_1}^t |J_1+J_2+J_3|(s_1,s_2)\d s_2\le \ff c t,\ \ T\ge t\ge 1, m\ge 1.$$
Combining this with  \eqref{007} and  \eqref{EG}, we find constants $c_3,c_4,c_5,c_6>0$ such that 
\beg{align*} & \sup_{T\ge t}  \Big| t  \E^\nu[\big[\mu_0(|\nn L_0^{-1} (\rr_{t,r}-1)|^2)\big|T<\tau\big]- \sum_{m=0}^\infty \ff{\e^{-(\ll_m-\ll_0) r }}{(\ll_m-\ll_0)^2}\Big|\\
&\le \ff{c_3}  t \sum_{m=1}^\infty \ff{\e^{-2(\ll_m-\ll_0)r} }{\ll_m-\ll_0} \le \ff {c_4}t \int_1^\infty s^{-\ff 2 d} \e^{-c_5 s^{\ff 2 d} r} \d s \le c_6 t^{-1}  \big( r^{-\ff{(d-2)^+}2} +1_{\{d=2\}} \log r^{-1}\big),\ \ t\ge 1.\end{align*}
 \end{proof} 
 
 \beg{lem}\label{LEM} There  exists a constant $c>0$ such that for any $t>0$ and nonnegative random variable $\xi\in \si(X_s:s\le t)$,
 $$\sup_{T\ge t} \E^\nu[\xi |T<\tau]\le c \E^\nu[\xi  |t<\tau],\ \ t\ge 1, \nu\in \scr P_0.$$\end{lem}
 \beg{proof} 
By the Markov property,  \eqref{AC0} for $p=q=\infty$ and \eqref{D10}, we find constants $c_1,c_2>0$ such that 
\beg{align*} &\E^\nu [\xi 1_{\{T<\tau\}}]= \E^\nu[ \xi 1_{\{t<\tau\}} P_{T-t}^D1(X_t)]\le c_1\e^{-\ll_0(T-t)} \E^\nu[\xi1_{\{t<\tau\}}],\\
&\P^\nu(T<\tau) \ge c_2 \P^\nu(t<\tau) \e^{-(T-t)\ll_0},\ \ T\ge t\ge 1.\end{align*} 
 Then 
 $$\E^\nu[\xi |T<\tau]=\ff{\E^\nu [\xi 1_{\{T<\tau\}}]}{\P^\nu(T<\tau} \le\ff{c_1 \E^\nu[\xi1_{\{t<\tau\}}]}{c_2\P^\nu(t<\tau)}=\ff{c_1}{c_2} \E^\nu[\xi|t<\tau].$$
 \end{proof}
 
 \beg{lem}\label{LL2}  Let $d\le 3$ and   denote $\nu_0=\ff{\phi_0}{\mu(\phi_0)}\mu.$   For any $\vv\in (\ff d 4\lor \ff{d^2}{2d +4}, 1)\ne\emptyset$, there exists a 
  constant  $c>0$  such that  
  $$   \sup_{T\ge t} \E^{\nu_0}\big[|\rr_{t,r}(y)-1|^{2}\big|T<\tau\big]\le c \phi_0^{-2}  (y) t^{-1} r^{-\vv},\ \ t\ge 1, r\in (0,1], y\in M^{\circ}.$$\end{lem}
  
  \beg{proof}  By Lemma \ref{LEM}, it suffices to prove for $T=t$ replacing $T\ge t$. For fixed $y\in M^\circ$, let $f= p_r^0(\cdot,y)-1.$ We have
 $$\rr_{t,r}(y)-1= \ff 1 t \int_0^t f(X_s)\d s.$$
  Then 
\beq\label{SS03}  \E^{\nu_0}\big[|\rr_{t,r}(y)-1|^21_{\{t<\tau\}}\big]= \ff 2 {t^2} \int_0^t \d s_1 \int_{s_1}^t \E^{\nu_0}\big[1_{\{t<\tau\}} f(X_{s_1})f(X_{s_2})\big]\d s_2.\end{equation}
  By \eqref{SD2},   $\mu_0(f)=0$,  and   the symmetry of  $P_t^0$ in $L^2(\mu_0), $ we obtain
\beq\label{NDD}    \beg{split} & I:=\e^{\ll_0 t} \E^{\nu_0}\big[1_{\{t<\tau\}} f(X_{s_1})f(X_{s_2})\big]=\mu(\phi_0)^{-1}  \mu_0\big(P_{s_1}^0\{fP_{s_2-s_1}^0(fP_{t-s_2}^0\phi_0^{-1})\}\big)\\
  &=\mu(\phi_0)^{-1}  \mu_0\big( fP_{s_2-s_1}^0(fP_{t-s_2}^0\phi_0^{-1}) \big)  = \mu(\phi_0)^{-1}  \mu_0\big( \{fP_{t-s_2}^0\phi_0^{-1}\} P_{s_2-s_1}^0 f \big)\\
  &= \mu(\phi_0)^{-1}  \mu_0\big( \{fP_{t-s_2}^0\phi_0^{-1}\}  \{P_{s_2-s_1}^0-\mu_0\} f \big).  \end{split} \end{equation}
Taking $q\in (\ff 5 2, 3)$ so that $\vv_1:=\ff{d+2}{2q}<1$ for $d\le 3$ and $\|\phi_0^{-1}\|_{L^q(\mu_0)}<\infty$ due to \eqref{JSS},   for any $p\in (1,2]$ we deduce from this and  \eqref{**1}   that    
 \beq\label{ND0} \beg{split}  &\mu(\phi_0)I\le    \|f\|_{L^p(\mu_0)} \|P_{t-s_2}^0\phi_0^{-1}\|_{L^\infty(\mu_0)} \|(P_{s_2-s_1}^0-\mu_0)f \|_{L^{\ff p{p-1}}(\mu_0)}\\ 
&\le    \|f\|_{L^p(\mu_0)}  \|P_{t-s_2}^0\|_{L^q(\mu_0)\to L^\infty(\mu_0)}\|\phi_0^{-1}\|_{L^q(\mu_0)} \|P_{s_2-s_1}^0-\mu_0\|_{L^2(\mu_0)\to L^{\ff p {p-1}}(\mu_0)} \|f\|_{L^2(\mu_0)}\\
&\le c_1 \|f\|_{L^p(\mu_0)}\|f\|_{L^2(\mu_0)} \{1\land (t-s_2)\}^{-\vv_1} \{1\land (s_2-s_1)\}^{-\ff{(d+2)(2-p)}{2p}} \e^{-(\ll_1-\ll_0)(s_2-s_1)}\end{split} \end{equation}
 holds for some constants $c_1>0$. Since $f= p_r^0(\cdot,y)-1$ and $\inf \phi_0^{-1}>0$, by \eqref{AC0} and \eqref{PR0}, we find   constants $\bb_1,\bb_2>0$ such that  
 \beg{align*} &\|f\|_{L^p(\mu_0)}\le 1 +\|p_r^0(\cdot,y)\|_{L^p(\mu_0)} \le 1 + \e^{r\ll_0} \phi_0^{-1}(y) \|\phi_0^{-1} p_r^D(\cdot,y)\|_{L^p(\mu_0)}\\
 &\le 1+ \bb_1 \phi_0^{-1}(y) \|\phi_0\|_{\infty}^{\ff{2-p}p}\|p_r^D(\cdot,y)\|_{L^p(\mu)}\le \bb_2 \phi_0^{-1} (y) r^{-\ff{d(p-1)}{2p}},\ \ r\in (0,1],p\in [1,2].\end{align*}
Combining this with \eqref{ND0} we find a constant $c_2>0$ such that  
$$I\le c_2 \phi_0^{-2}(y)r^{-\ff{d(p-1)}{2p}- \ff{d}4 }  \{1\land (t-s_2)\}^{-\vv_1} \{1\land (s_2-s_1)\}^{-\ff{(d+2)(2-p)}{2p}} \e^{-(\ll_1-\ll_0)(s_2-s_1)},\ \ p\in (1,2].$$
Taking $p> p_0:=1\lor \ff{2(d+2)}{d+6}$ such that  
$$\ \ \vv_2:= \ff{(d+2)(2-p)}{4p}\le  \ff{5(2-p)}{4p}<1,$$
we arrive at
$$I\le c_2\phi^{-2}(y) r^{-\ff{d(p-1)}{2p}- \ff{d}4 }  \{1\land (t-s_2)\}^{-\vv_1} \{1\land (s_2-s_1)\}^{- \vv_2} \e^{-(\ll_1-\ll_0)(s_2-s_1)}$$ for some constants $\vv_1,\vv_2\in (0,1)$. 
Combining this with \eqref{SS03}, we obtain 
$$  \E^{\nu_0}\big[|\rr_{t,r}(y)-1|^{2}\big|t<\tau\big]\le c\phi_0^{-2}  (y) t^{-1} r^{-\ff{d(p-1)}{2p}- \ff{d}4},\ \  t\ge 1. $$
Noting that 
$$\lim_{p\downarrow p_0}  \Big\{\ff{d(p-1)}{2p}+\ff d 4 \Big\}= \ff d 4\lor \ff{d^2}{2d +4}<1\ \text{for}\ d\le 3,$$ for any $\vv\in (\ff d 4\lor \ff{d^2}{2d +4}, 1)$, there exists $p>p_0$ such that $\ff d 4\lor \ff{d^2}{2d +4}\le \vv$.
Therefore, the proof is finished. 
 \end{proof}

  \beg{lem}\label{LL2'}  Let $d\le 3$ and denote  $\psi_m(t)=\ff 1 t \int_0^t (\phi_m\phi_0^{-1})(X_s)\d s.$ Then
  there exists a constant $c>0$ such that  for any $p\in [1,2],$
  $$ \sup_{T\ge t}    \E^{\nu_0}\big[|\psi_m(t)|^{2p}\big|t<\tau\big]\le c m^{\ff{p(d+4)-d-8}{2d}} t^{-p},\ \ t\ge 1, m\ge 1, r\in (0,1).$$\end{lem}
  
  \beg{proof}   By Lemma \ref{LEM}, it suffices to prove for $T=t$ replacing $T\ge t$.  By H\"older's inequality, we have 
  $$\E^{\nu_0}\big[|\psi_m(t)|^{2p}|T<\tau\big]\le \big\{\E^{\nu_0}\big[|\psi_m(t)|^{2}|T<\tau\big]\big\}^{2-p} \big\{\E^{\nu_0}\big[|\psi_m(t)|^{4}|T<\tau\big]\big\}^{p-1}.$$
 Combining this with \eqref{D10}, 
   it suffices to find a constant $c>0$ such that 
   \beq\label{SS01'}    \E^{\nu_0}\big[|\psi_m(t)|^21_{\{t<\tau\}}\big]  \le \ff{c \e^{-\ll _0t} }{tm^{\ff 2 d}} ,\ \  t\ge 1, r\in (0,1),\end{equation}
  \beq\label{SS02'}  \E^{\nu_0}\big[|\psi_m(t)|^41_{\{t<\tau\}}\big] \le c\ss m \, \e^{-\ll_0 t} t^{-2},\ \   t\ge 1, r\in (0,1).\end{equation}

 (a) Proof of \eqref{SS01'}.   Let $\hat \phi_m= \phi_m\phi_0^{-1}.$ We have 
\beq\label{SS03'}  \E^{\nu_0}\big[|\psi_m(t)|^21_{\{t<\tau\}}\big]= \ff 2 {t^2} \int_0^t \d s_1 \int_{s_1}^t \E^{\nu_0}\big[1_{\{t<\tau\}} \hat \phi_m (X_{s_1})\hat \phi_m(X_{s_2})\big]\d s_2.\end{equation}
  By \eqref{PR1},  \eqref{SD2},  $\mu_0(|\hat \phi_m|^2)=1$,    and the symmetry of $P_t^0$ in $L^2(\mu_0)$, we find a constant $c_1 >0$ such that 
  \beg{align*} &\e^{\ll_0 t} \E^{\nu_0}\big[1_{\{T<\tau\}} \hat \phi_m(X_{s_1})\hat \phi_m(X_{s_2})\big]= \nu_0\big(\phi_0 P_{s_1}^0\{\hat \phi_m P_{s_2-s_1}^0(\hat \phi_m P_{t-s_2}^0\phi_0^{-1})\}\big)\\
  &=\ff 1 {\mu(\phi_0)}  \mu_0\big(\hat \phi_m P_{s_2-s_1}^0(\hat \phi_m P_{t-s_2}^0 \phi_0^{-1})\big)=\ff{ \e^{-(\ll_m-\ll_0)(s_2-s_1)} }{\mu(\phi_0)} \mu_0\big(|\hat \phi_m|^2    P_{t-s_2}^0 \phi_0^{-1})\big)\\
 & \le c_1 \e^{-(\ll_m-\ll_0)(s_2-s_1)}   \|P_{t-s_2}\|_{L^p(\mu_0)\to\infty(\mu_0)}\|\phi_0^{-1}\|_{L^p(\mu_0)},\ \ p>1.\end{align*}
  Since $d\le 3$, we may take $p\in (1,3)$ such that $\vv:=\ff{d+2}{2q}<1$ and $\|\phi_0^{-1}\|_{L^p(\mu_0)}<\infty$ due to \eqref{JSS}, so that this and \eqref{**1} imply  
  $$\e^{\ll_0 t} \E^{\nu_0}\big[1_{\{t<\tau\}} \hat \phi_m(X_{s_1})\hat \phi_m(X_{s_2})\big]\le  c_2 \e^{-(\ll_m-\ll_0) (s_2-s_1)} \{1\land (t-s_2)\}^{-\vv} $$
  for some constant $c_3>0$. Therefore, \eqref{SS01'} follows from \eqref{SS03'} and \eqref{EG}. 
  
   (b) Proof of \eqref{SS02'}.  For any $s>0$ we have 
 \beq\label{SS04'} \beg{split} &s^4 \E^{\nu_0}\big[|\psi_m(s)|^{4}1_{\{s<\tau\}}\big]\\
& = 24  \int_0^s \d s_2\int_{s_1}^s \d s_2 \int_{s_2}^s \d s_3 \int_{s_3}^s\E^{\nu_0}\big[1_{\{s<\tau\}} \hat\phi_m(X_{s_1})\hat\phi_m(X_{s_2})\hat\phi_m(X_{s_3}) \hat\phi_m(X_{s_4}) \big]\d s_4\\
 &= 24  \int_0^s \d s_2\int_{s_1}^s \d s_2 \int_{s_2}^s \d s_3 \int_{s_3}^s\E^{\nu_0}\big[1_{\{s_3<\tau\}} \hat\phi_m(X_{s_1})\hat\phi_m(X_{s_2})g_s(s_3,s_4)\big]\d s_4,\end{split}\end{equation}
 where due to \eqref{SD2} and the Markov property,    
 \beq\label{SS05'} \beg{split} g_s(s_3,s_4)&:= \E^{\nu_0}\big[1_{\{s<\tau\}} \hat\phi_m(X_{s_3})\hat\phi_m(X_{s_4})\big|X_r: r\le s_3\big]\\  
 &= \hat\phi_m(X_{s_3}) \E^{X_{s_3}} \big[1_{\{s-s_3<\tau\}}\hat\phi_m(X_{s_4-s_3})\big]\\ 
 &=\e^{-\ll_0 (s-s_3)} \big\{ \hat\phi_m \phi_0  P_{s_4-s_3}^0 (\hat\phi_mP_{s-s_4}^0 \phi_0^{-1})\big\}(X_{s_3}),\ \ 0<s_3<s_4\le s.\end{split}\end{equation} 
 So, by Fubini's theorem and Schwarz's inequality,  we obtain 
 \beg{align*} &I(s):= s^4 \e^{\ll_0s}  \E^{\nu_0}\big[|\psi_m(s)|^{4}1_{\{s<\tau\}}\big]\\
 &= 12 \e^{\ll_0 s} \int_0^s \d r_1 \int_{r_1}^s \E^{\nu_0}\bigg[ 1_{\{r_1<\tau\}}  g_s(r_1,r_2) \bigg|\int_0^{r_1} \hat\phi_m(X_r)\d r\bigg|^2 \bigg]\d r_2\\
 &\le 12 \sup_{r\in [0,s]} \ss{I(r)} \int_0^s \d r_1 \int_{r_1}^s \Big\{\e^{2\ll_0s -\ll_0r_1} \E^{\nu_0}\big[1_{\{r_1<\tau\}}g_s(r_1,r_2)^2\big]\Big\}^{\ff 1 2 } \d r_2.\end{align*}
Consequently, 
\beq\label{QK0'} I(t)\le \sup_{s\in [0,t]}I(s)\le \bigg(12  \sup_{s\in [0,t]} \int_0^s \d r_1 \int_{r_1}^s \Big\{\e^{\ll_0(2s - r_1)} \E^{\nu_0}\big[1_{\{r_1<\tau\}}g_s(r_1,r_2)^2\big]\Big\}^{\ff 1 2 } \d r_2\bigg)^2.\end{equation} 

On the other hand,   by the definition of $\nu_0$,   \eqref{SD2},   \eqref{SS05'}  and that $\mu_0$ is $P_t^0$-invariant,  we obtain 
\beq\label{L009}\beg{split}&\E^{\nu_0}\big[1_{\{r_1<\tau\}}|g_s(r_1,r_2)|^2\big]\\ 
&\le \ff{\e^{-2\ll_0(s-r_1)-\ll_0r_1}}{\mu(\phi_0)}   \mu_0\big(P_{r_1}^0\{\phi_0^{-1}| \hat\phi_m \phi_0  P_{r_2-r_1}^0 (\hat\phi_mP_{s-r_2}^0 \phi_0^{-1})|^2\}\big)\\  
&=  \ff{\e^{-\ll_0(2s-r_1) }}{\mu(\phi_0)} \mu_0\big(\phi_0| \hat\phi_m   P_{r_2-r_1}^0 (\hat\phi_mP_{s-r_2}^0 \phi_0^{-1})|^2\big)\\   
&\le    \ff{2\e^{-\ll_0(2s-r_1)}}{\mu(\phi_0)}\mu_0\big(\phi_0\{| \hat\phi_m   (P_{r_2-r_1}^0 \hat\phi_m) \mu(\phi_0)|^2 + | \hat\phi_m   P_{r_2-r_1}^0(\hat\phi_m[P_{s-r_2}^0-\mu_0]  \phi_0^{-1})|^2\}\big).\end{split}\end{equation} 
Then,  by \eqref{SS05'},  \eqref{PR1},  \eqref{SD2},  $\mu_0(|\hat \phi_m|^2)=1$,    and  noting that $\mu_0$ is $P_t^0$-invariant, we find a constant $c_1 >0$ such that
\beg{align*}&\E^{\nu_0}\big[1_{\{r_1<\tau\}}|g_s(r_1,r_2)|^2\big]\le 2  \e^{-\ll_0(2s-r_1)-(\ll_m-\ll_0)(r_2-r_1)} \|\phi_m\|_\infty\|\phi_0\|_\infty \mu_0(|\hat \phi_m| |P_{(r_2-r_1)/2}\hat \phi_m|^2)  \\     
&\qquad +  2  \ff{\e^{-\ll_0(2s-r_1)}\|\phi_m\|_\infty}{\mu(\phi_0)}
 \mu_0\big( |\hat\phi_m| \cdot\big|  P_{r_2-r_1}^0 (\hat\phi_m(P_{s-r_2}^0-\mu_0) \phi_0^{-1})\big|^2\big)   \\  
 &\le c_1 \e^{-\ll_0(2s-r_1)}\Big\{  \e^{-(\ll_m-\ll_0)(r_2-r_1)}\|\phi_m\|_\infty \|P_{(r_2-r_1)/2}-\mu_0\|_{L^2(\mu_0)\to L^4(\mu_0)}^2  \\
 &\qquad \qquad +  \|\phi_m\|_\infty  \|P_{r_2-r_1}^0(\hat\phi_m[P_{s-r_2}^0-\mu_0]  \phi_0^{-1})\|^2_{L^4(\mu_0)}\Big\}.\end{align*} 
By \eqref{EG},   \eqref{**1},   $\|\hat\phi_m\|_{L^2(\mu_0)}=1$,  $\|\phi_0^{-1}\|_{L^q(\mu_0)}<\infty$  and  $\vv:=\ff{d+2}8\lor \ff{d+2}{2q}<1$ for $q\in (\ff 5 2,3)$ due to \eqref{JSS} and  $d\le 3$,  we find constants $c_2>0$   such that   
  $$  \|\phi_m\|_\infty \|P_{(r_2-r_1)/2}-\mu_0\|_{L^2(\mu_0)\to L^2(\mu_0)}^2\le c_2 \ss m \{1\land (r_2-r_2)\}^{-\ff d 4},$$ and 
\beg{align*}& \|\phi_m\|_\infty  \|P_{r_2-r_1}^0(\hat\phi_m[P_{s-r_2}^0-\mu_0]  \phi_0^{-1})\|^2_{L^4(\mu_0)}\\  
 & \le \|\phi_m\|_\infty\|P_{r_2-r_1}^0  \|_{L^2(\mu_0)\to L^4(\mu_0)}^2 \|\hat\phi_m\|_{L^2(\mu_0)}^2\|(P_{s-r_2}^0-\mu_0) \phi_0^{-1}\|_{L^\infty(\mu_0)}^2   \\ 
 & \le \|\phi_m\|_\infty\|P_{r_2-r_1}^0  \|_{L^2(\mu_0)\to L^4(\mu_0)}^2 \| P_{s-r_2}^0-\mu_0\|_{L^q(\mu_0)\to L^\infty(\mu_0)}^2\|\phi_0^{-1}\|_{L^q(\mu_0)}^2\\   
  &\le c_2 \ss m\, \e^{-2(\ll_1-\ll_0)(s-r_2)} \{1\land (r_2-r_1)\}^{-2\vv} \{1\land (s-r_2)\}^{-2\vv}.\end{align*}
 Therefore, there exist  constants $c_3>0$ and $\vv\in (0,1)$ such that 
\beg{align*} & \E^{\nu_0}\big[1_{\{r_1<\tau\}}|g_s(r_1,r_2)|^2\big]\le  c_3 \e^{-\ll_0(2s-r_1) -(\ll_m-\ll_0)(r_2-r_1)}\ss m \{1\land (r_2-r_2)\}^{-\ff d 4}\\
 &\qquad\qquad +  c_3\ss m\,\e^{-\ll_0(2s-r_1)-2(\ll_1-\ll_0)(s-r_2)}  \{1\land (r_2-r_1)\}^{-2\vv} \{1\land (t-r_2)\}^{-2\vv}.\end{align*} 
 Combining this with \eqref{QK0'} and the definition of $I(t)$,  we prove \eqref{SS02'} for some constant $c>0,$ and hence finish the proof. 
   \end{proof} 
   
\beg{lem}\label{LL4} Let $d\le 3$. Then for any $p\in (1,  \ff{3d+16}{5d+8}\land \ff{d+2}{d+1})\ne\emptyset$,   there exists a constant $c>0$ such that
$$ \sup_{r>0,T\ge t} \E^{\nu_0}\big[ \mu_0( |\nn L_0^{-1} (\rr_{t,r}-1)|^{2p})|T<\tau\big]\le c t^{-p},\ \ t\ge 1.$$
\end{lem} 

\beg{proof}   By Lemma \ref{LEM}, it suffices to prove for $T=t$ replacing $T\ge t$. Let  $p\in (1,  \ff{3d+16}{5d+8}\land \ff{d+2}{d+1})$, where $p>1$ is equivalent to
\beq\label{PPP} \ff{p}{2p-1}<1,  \end{equation}  while $p< \ff{3d+16}{5d+8}\land \ff{d+2}{d+1}$ implies 
$$\ff{(d+2)(2p-2)}4+\ff{d(p-1)} 2 +\Big(\ff{p(d+4)+d} 4 -2\Big)^+<1,$$ and hence 
  there exists $\vv\in (0,1)$ such that 
\beq\label{PPP'} \ff{(d+2)(2p-2+\vv)} 4+ \ff{d(p-1)} 2 +\Big(\ff{p(d+4)+d} 4 -2\Big)^+<1.\end{equation} 
By \eqref{**1},   \eqref{GRD}, $L_0^{-1}= -\int_0^\infty P_s^0\d s,$  and applying H\"older's inequality, we find a constant $c_1,c_2>0$ such that 
 \beq\label{MFN} \beg{split} &  \int_M \big|\nn L_0^{-1} (\rr_{t,r}-1)\big|^{2p}\d\mu_0 \le \int_M \bigg(\int_0^\infty \big|\nn P_s^0 (\rr_{t,r}-1)\big|\d s \bigg)^{2p}\d\mu_0\\
&\le c_1 \int_M \bigg(\int_0^\infty \ff 1 {\ss s} \big\{P_{\ff s 4}^0 \big|P_{\ff {3s} 4}^0 (\rr_{t,r}-1) \big|^p\big\}^{\ff 1 p} \d s   \bigg)^{2p}\phi_0^{-\vv} \d\mu_0\\
&\le c_1 \bigg(\int_0^\infty s^{-\ff{p}{2p-1}} \e^{-\ff{2p\theta s}{2p-1}}\d s\bigg)^{\ff {2p-1}{2p}} \int_0^\infty \e^{\theta s} \mu_0\big(\phi_0^{-\vv} \big\{P_{\ff s 4}^0 |P_{\ff {3s} 4}^0 (\rr_{t,r}-1)|^p\big\}^2\big)\d s,\ \  \theta>0.\end{split} \end{equation}
Noting that   $\ff{p}{2p-1}<1$ due to \eqref{PPP}, we obtain
\beq\label{YP} \int_0^\infty s^{-\ff{p}{2p-1}} \e^{-\ff{2p\theta s}{2p-1}}\d <\infty,\ \ \theta>0.\end{equation}
Moreover, since $\|\phi_0^{-\vv}\|_{L^{2\vv^{-1}}(\mu_0)}=1, \mu_0(\rr_{t,r}-1)=0$, and $P_t^0$ is contractive in $L^p(\mu_0)$ for $p\ge 1$,     by \eqref{**1} and H\"older's inequality, we find a constant $c_2>0$ such that 
\beg{align*} & \mu_0\big(\phi_0^{-\vv} \big\{P_{\ff s 4}^0 |P_{\ff {3s} 4}^0 (\rr_{t,r}-1)|^p\big\}^2\big)\le \big\|P_{\ff s 4}^0|P_{\ff {3s} 4}^0 (\rr_{t,r}-1)|^p\big\|_{L^{\ff 4 {2-\vv}}(\mu_0)}^2\|\phi_0^{-\vv}\|_{L^{2\vv^{-1}}(\mu_0)}\\
&\le \|P_{\ff s 4}^0\|_{L^{\ff 4{2-\vv}}(\mu_0)}^2\big\|(P_{\ff s 2}^0-\mu_0)(P_{\ff {s} 4}^0 \rr_{t,r}-1)|\big\|_{L^{\ff{4p}{2-\vv}}(\mu_0)}^{2p} \\
&\le \|P_{\ff s 2}^0-\mu_0\|_{L^2(\mu_0)\to L^{\ff{4p}{2-\vv}}(\mu_0)}^{2p} \|P_{\ff s 4}^0\rr_{t,r}-1\|_{L^2(\mu_0)}^{2p}\\
&\le c_2 (1\land s)^{-\ff{(d+2)(2p-2+\vv)}{4}}\e^{-(\ll_1-\ll_0)ps} \|P_{\ff s 4}^0\rr_{t,r}-1\|_{L^2(\mu_0)}^{2p}.\end{align*}
Combining this with \eqref{YP}, we find a function $c: (0,\infty)\to (0,\infty)$  such that 
\beq\label{MNN}  \beg{split} &\E^{\nu_0}\big[1_{\{t<\tau\}}\mu_0(  |\nn L_0^{-1} (\rr_{t,r}-1) |^{2p} ) \big] \\
&\le c(\theta) \int_0^\infty \e^{\theta s} (1\land s)^{-\ff{(d+2)(2p-2+\vv)} 4} \e^{-(\ll_1-\ll_0)ps} \E^{\nu_0}\big[1_{\{t<\tau\}} \|P_{\ff s 4}^0\rr_{t,r}-1\|_{L^2(\mu_0)}^{2p}\big]\d s,\ \ \theta>0.\end{split}\end{equation} 
  By \eqref{PR1}, \eqref{DS1} and H\"older's inequality, we obtain
\beg{align*} & \|P_{\ff s 4}^0\rr_{t,r}-1\|_{L^2(\mu_0)}^{2p} = \Big(\sum_{m=1}^\infty \e^{-(\ll_m-\ll_0)(2r+s/2)}|\psi_m(t)|^2\Big)^{p}\\
&\le \Big(\sum_{m=1}^\infty \e^{-(\ll_m-\ll_0)(2r+s/2)}\Big)^{p-1} \sum_{m=1}^\infty  \e^{-(\ll_m-\ll_0)(2r+s/2)}|\psi_m(t)|^{2p}.\end{align*}
Noting that \eqref{EG} implies
$$\sum_{m=1}^\infty \e^{-(\ll_m-\ll_0)(2r+s/2)}\le a_1 \int_1^\infty \e^{-\aa_2(r+s/2)t^{\ff 2 d}}\d t\le \aa_3 (1\land s)^{-\ff d 2} $$
for some constants $\aa_1,\aa_2,\aa_3>0$, we derive 
 $$ \E^{\nu_0}\big[  \|P_{\ff s 4}^0\rr_{t,r}-1\|_{L^2(\mu_0)}^{2p} \big|t<\tau \big] \le c_3(1\land s)^{-\ff{d(p-1)} 2}    \sum_{m=1}^\infty \e^{-(\ll_m-\ll_0)(2r+s/2)} \E^{\nu_0}\big[|\psi_m(t)|^{2p}\big|t<\tau\big]$$
 for some constant $c_3>0$.  Combining this with Lemma \ref{LL2'},  \eqref{EG}, we find constants $c_4,c_5,c_6,c_7>0$ such that
 \beg{align*} & \E^{\nu_0}\big[ \|P_{\ff s 4}^0\rr_{t,r}-1\|_{L^2(\mu_0)}^{2p}\big|t<\tau \big] \le c_4t^{-p} (1\land s)^{-\ff{d(p-1)}2}  \int_1^\infty   
 \e^{-c_5 su^{\ff 2 d}}  u^{\ff{p(d+4)-d-8}{2d}}\d u\\
 &\le  c_6 t^{-p} (1\land s)^{-\ff{d(p-1)}2} s^{2-\ff{p(d+4)+d}4}\int_s^\infty t^{\ff{p(d+4)+d}4-3} \e^{-t} \d t\\
 &\le c_7 t^{-p}(1\land s)^{-\ff{d(p-1)}2- (\ff{p(d+4)+d}4-2)^+}  \log (2+s^{-1}),\end{align*}  where the term $\log (2+s^{-1})$ comes when $\ff{p(d+4)+d}4-3=-1$. 
This together with \eqref{PPP'} and \eqref{MNN} for  $\theta \in (0,  \ll_1-\ll_0)$ implies the desired estimate. 
 \end{proof}

\beg{lem}\label{LL3} Let $d\le 3$. If $r_t= t^{-\aa}$ for some $\aa\in (1,\ff 4 d\land \ff{2d+4}{d^2})\ne\emptyset$, then  $\rr_{t,r_t,r_t}:= (1-r_t)\rr_{t,r_t} + r_t$ satisfies 
$$\lim_{t\to\infty}   \sup_{T\ge t} \E^{\nu_0}\big[\mu_0(|\scr M(\rr_{t,r_t,r_t}, 1)^{-1}-1|^{q})\big|T<\tau\big] =0,\ \ q\ge 1.$$
\end{lem} 

\beg{proof} By Lemma \ref{LEM}, it suffices to prove for $T=t$ replacing $T\ge t$.  By the same reason leading to (3.16) in \cite{WZ20}, for any $\eta\in (0,1), y\in M$, we have  
 $$  \E^{\nu_0}\big[|\scr M(\rr_{t,r_t,r_t}(y), 1)^{-1}-1|^{q}\big|t<\tau\big]\le \Big|\ff 1 {\ss{1-\eta}} - \ff 2 {2+\eta}\Big|^{q} + \P^{\nu_0}\big(|\rr_{t,r_t}(y)-1|>\eta\big).$$ 
Combining this with Lemma \ref{LL2}  we find   constants $c>0$ and $\vv \in (0, \aa^{-1})$ such that 
$$\E^{\nu_0}\big[|\scr M(\rr_{t,r_t,r_t}(y), 1)^{-1}-1|^{q}\big|t<\tau\big]\le \Big|\ff 1 {\ss{1-\eta}} - \ff 2 {2+\eta}\Big|^{q} +  c \eta^{-1} \phi_0(y)^{-2} t^{-1+\aa \vv }.$$
Since  $\mu_0(\phi_0^{-2})=1$, 
we obtain
$$ \E^{\nu_0}\big[\mu_0(|\scr M(\rr_{t,r_t,r_t}, 1)^{-1}-1|^{q})\big|t<\tau\big] \le  \Big|\ff 1 {\ss{1-\eta}} - \ff 2 {2+\eta}\Big|^{q} +  c \eta^{-1} t^{-1+\aa \vv },\ \ \eta\in (0,1), t\ge 1.$$
Noting that $\aa\vv<1$,  by letting first $t\to \infty$ then $\eta\to 0$,  we finish the proof. 
\end{proof}

\beg{lem}\label{LL5} Let $\mu_{t,r,r}= (1+\rr_{t,r,r})\mu_0$, where $\rr_{t,r,r}:= (1-r) \rr_{t,r}+r, r\in (0,1].$ Assume that $\nu= h\mu$  with $h\phi_0^{-1} \in L^p(\mu_0)$ for some $p>1$.  Then there exists a constant $c>0$ such that
$$ \sup_{T\ge t} \E^\nu\big[\W_2(\mu_{t,r,r},\mu_{t})^2\big|T<\tau\big] \le c r,\ \ t>0, r\in (0,1].$$
\end{lem} 
\beg{proof} By Lemma \ref{LEM}, it suffices to prove for $T=t$ replacing $T\ge t$.  Firstly,   it is easy to see that 
\beq\label{*N0} \W_2(\mu_{t,r,r}, \mu_{t,r})^2\le D^2\|\mu_{t,r,r}- \mu_{t,r}\|_{var}=D^2 \mu_0(|\rr_{t,r,r}-\rr_{t,r}|)\le 2D^2r,\ \ r\in (0,1].\end{equation}
Next, by the definition of $\mu_{t,r}$,  we have 
$$\pi(\d x,\d y):= \mu_t(\d x) P_r^0(x,\d y)\in \C(\mu_t,\mu_{t,r}),$$
where $P_r^0(x,\cdot)$ is the distribution of $X_r^0$ starting at $x$. So, 
\beq\label{*N2} \W_2(\mu_t,\mu_{t,r})^2\le \int_M \E^x[ \rr(x,X_r^0)^2 ] \mu_t(\d x).\end{equation}
Moreover, by It\^o's formula and $L_0=L+2\nn\log \phi_0$, we find a constant $c_1>0$ such that 
$$\d \rr(x, X_r^0)^2 = L_0 \rr(x,\cdot)^2(X_r^0) \d r + \d M_r\le\big\{ c_1 + c_1 \phi_0^{-1} (X_r^0)\big\}\d r+\d M_r$$
holds for some martingale $M_r$.  Combining this with \eqref{*B*}, and noting that $\log(1+\phi_0^{-1})\ge \log (1+\|\phi_0\|_\infty^{-1})>0$, we find a constant $c_2>0$ such that 
\beg{align*}&\W_2(\mu_t,\mu_{t,r})^2\le c_1 r + c_1 \int_M \bigg(\E^x \int_0^r  \phi_0^{-1} (X_s^0)\d s \bigg)\mu_t(\d x) \\
&\le  c_2 r\mu_t(\log (1+\phi_0^{-1}))= \ff {c_2r} {t} \int_0^t \log \{1+\phi_0^{-1} (X_s)\}\d s, \ \ r\in (0,1].\end{align*} 
Combining this with \eqref{*N0}, \eqref{SD2},  $\|P_t^0\|_{L^p(\mu_0)}=1$ for $t\ge 0$ and $p\ge 1$, and noting that 
$$ \inf_{t\ge 0}  \mu_0(h\phi_0^{-1} P_t^0\phi_0^{-1}) >0,$$ we find  constants $c_3,c_4>0$ such that 
\beq\label{US} \beg{split} &\E^\nu[\W_2(\mu_{t,r,r},\mu_t)^2|t<\tau] =\ff{\E^\nu[1_{\{t<\tau\}}\W_2(\mu_{t,r,r},\mu_t)^2]}{\P^\nu(t<\tau)}\\
&\le   \ff {c_3r}{t\mu_0(h\phi_0^{-1}P_t^0\phi_0^{-1})}\int_0^t  \mu_0(h\phi_0^{-1}  P_s^0\log\{1+\phi_0^{-1}\})\d s\\
&\le c_3 r \|h\phi_0^{-1}\|_{L^p(\mu_0)}\|\log(1+\phi_0^{-1})\|_{L^{\ff p{p-1}}(\mu_0)}\le c_4 r,\ \ r\in (0,1].\end{split}\end{equation}
Combining this with \eqref{*N0} we finish the proof. 

\end{proof} 

We are now ready to prove the main result in this section. 

\beg{proof}[Proof of Proposition $\ref{Pn1}(1)$]  Since the upper bound is infinite for $d\ge 4$, it suffices to consider $d\le 3.$   

(a)   We first assume that $\nu= h\mu$  with $h\le C\phi_0$ for some constant $C>0$. In this case, by \eqref{D10} and $\E^\nu=\int_M\E^x\nu(\d x)$, there exists a constant $c_0>0$ such that
\beq\label{EQ0} \E^\nu(\cdot|t<\tau)\le c_0 \E^{\nu_0}(\cdot|t<\tau),\ \ t\ge 1.\end{equation}
Let $\mu_{t,r_t,r_t}= \{(1-r_t) \rr_{t,r_t} + r_t\}\mu_0$ with  $r_t=t^{-\aa}$ for some $\aa\in  (1,\ff 4 d\land \ff{2d+4}{d^2})$. By Lemma \ref{LL5} and the triangle inequality of $\W_2$,  there exists a constant $c_1>0$ such that for any $t\ge 1$, 
\beq\label{LM1} \E^\nu\big[ \W_2(\mu_t,\mu_0)^2 \big|t<\tau\big]  \le (1+\vv)\E^\nu\big[ \W_2(\mu_{t,r_t,r_t}, \mu_0)^2\big|t<\tau\big] + c_1(1+\vv^{-1}) t^{-\aa},\ \ \vv>0.  \end{equation}
On the other hand, by \eqref{*UPA}, \eqref{EQ0}, Lemmas \ref{LL1}, \ref{LL4} and \ref{LL3},  there exists $p>1$ such that 
\beg{align*} &\limsup_{t\to\infty} t \E^\nu\big[ \W_2(\mu_{t,r_t,r_t},\mu_0)^2\big|t<\tau\big]\le \limsup_{t\to\infty} t\E^\nu\bigg[ \int_M \ff{|\nn L_0^{-1}(\rr_{t,r_t}-1)|^2}{\scr M(\rr_{t,r_t,r_t},1)} \d\mu_0 \bigg|t<\tau\bigg]\\
&\le \limsup_{t\to\infty} t \Big\{ \E^\nu\big[ \mu_0(|\nn L_0^{-1}(\rr_{t,r_t}-1)|^2) \d\mu_0 \big|t<\tau\big]\\
&\qquad + \big(\E^\nu\big[ \mu_0(|\nn L_0^{-1}(\rr_{t,r_t}-1)|^{2p}) \d\mu_0 \big|t<\tau\big]\big)^{\ff 1 p}\big( \E^{\nu}\big[\mu_0(|\scr M(\rr_{t,r_t,r_t}, 1)^{-1}-1|^{\ff p{p-1}})\big|t<\tau\big]\big)^{\ff {p-1}p}  \Big\}\\
&= \limsup_{t\to\infty}  t\E^\nu\big[ \mu_0(|\nn L_0^{-1}(\rr_{t,r_t}-1)|^2) \d\mu_0 \big|t<\tau\big]\le  \sum_{m=1}^\infty \ff 2 {(\ll_m-\ll_0)^2}.\end{align*}
Combining this with \eqref{LM1} where $\aa>1$, we prove  \eqref{R0}.   

(b) In general, for any $t\ge 2$ and $\vv\in (0,1)$, we consider 
$$\mu_t^\vv:= \ff 1 {t-\vv} \int_\vv^t \dd_{X_s}\d s.$$
Letting $D$ be the diameter of $D$, we find a constant $c_1>0$ such that 
\beq\label{NL1} \W_2(\mu_t^\vv,\mu_t)^2 \le D^2 \|\mu_t-\mu_t^\vv\|_{var} \le  c_1 \vv t^{-1},\ \ t\ge 2, \vv\in (0,1).\end{equation} 
On the other hand, by the Markov property we obtain
\beg{align*} &\E^\nu\big[1_{\{t<\tau\}}\W_2(\mu_t^\vv,\mu_0)^2\big]= \E^\nu\big[1_{\{\vv<\tau\}}  \E^{X_\vv} (1_{\{t-\vv<\tau\}} \W_2(\mu_{t-\vv},\mu_0)^2)\big]\\
&= \P^\nu(\vv<\tau) \E^{\nu_\vv} \big[1_{\{t-\vv<\tau\}} \W_2(\mu_{t-\vv},\mu_0)^2\big]\\
&=\P^{\nu_\vv}(t-\vv<\tau)\P^\nu(\vv<\tau) \E^{\nu_\vv} \big[\W_2(\mu_{t-\vv},\mu_0)^2\big|t-\vv<\tau\big],\end{align*}
where $\nu_\vv=h_\vv\mu$ with
$$h_\vv(y):= \ff 1 {\P^\nu(\vv<\tau)} \int_M  p_\vv^D(x,y)\nu(\d x) \le c(\vv,\nu)\phi_0(y)$$ 
for some constant $c(\vv,\nu)>0$. Moreover, by \eqref{D09}, \eqref{D10}  and $\nu_\vv=h_\vv\mu$, we have 
$$\lim_{t\to\infty}   \ff{\P^{\nu_\vv}(t-\vv<\tau)\P^\nu(\vv<\tau)}{\P^\nu(t<\tau)}=1.$$ 
 So, (a) implies 
\beg{align*} &\limsup_{t\to\infty}\Big\{ t  \E^\nu\big[ \W_2(\mu_t^\vv,\mu_0)^2\big|t<\tau\big] \Big\}\\
&= \limsup_{t\to\infty} \ff{\P^{\nu_\vv}(t-\vv<\tau)\P^\nu(\vv<\tau)}{\P^\nu(t<\tau)}\Big\{ t  \E^{\nu_\vv}\big[ \W_2(\mu_{t-\vv},\mu_0)^2\big|t-\vv <\tau\big]\Big\}\\
&\le \sum_{m=1}^\infty \ff 2 {(\ll_m-\ll_0)^2}. \end{align*} 
Combining this with \eqref{NL1},  we arrive at 
\beg{align*} &\limsup_{t\to\infty}\Big\{ t  \E^\nu\big[ \W_2(\mu_t,\mu_0)^2\big|t<\tau\big] \Big\}\\
& \le (1+\vv^{\ff 1 2}) \limsup_{t\to\infty}\Big\{ t  \E^{\nu}\big[ \W_2(\mu_t^\vv,\mu_0)^2\big|t<\tau\big] \Big\} + c_1\vv (1+\vv^{-\ff 1 2}) \\
&\le  (1+\vv^{\ff 1 2})\sum_{m=1}^\infty \ff 2 {(\ll_m-\ll_0)^2}+ c_1\vv (1+\vv^{-\ff 1 2}),\ \ \vv\in (0,1). \end{align*} 
By letting $\vv\to 0$, we derive \eqref{R0}. 
\end{proof} 

\beg{proof}[Proof of Proposition $\ref{Pn1}(2)$-$(3)$]  Let $d\ge 4$. By \eqref{NL1}, it suffices to prove the desired estimates for $\mu_t^1$ replacing $\mu_t$. Therefore, 
we may and do assume $\nu=h\mu$ with $\|h\phi_0^{-1}\|_\infty<\infty$. 
 Since 
$$\lim_{p\downarrow p_0} \Big\{\ff d 2 +\ff{(d+2)(p-1)}{2p}-2\Big\} =\ff{2(d-4)} 3,$$   
by Lemma \ref{LL1}(1), for any $k>\ff{2(d-4)} 3$, there exist constants $c_1,c_2>0$ such that 
\beg{align*} & t\E^\nu\big[\mu_0(|\nn L_0^{-1} (\rr_{t,r}-1)|^2\big|T<\tau\big]\le c_1 \sum_{m=1}^\infty \ff{\e^{-2(\ll_m-\ll_0)r}}{(\ll_m-\ll_0)^2} +c_1 t^{-1} r^{-k}\\
&\le c_2 \big\{1 + 1_{\{d=4\}} \log r^{-1} + t^{-1} r^{-k}\big\},\ \ r\in (0,1), t\ge 1, T\ge t.\end{align*}   
Combining this with the following inequality  due to  \cite[Theorem 2]{Ledoux} for $p=2$: 
$$\W_2(f\mu_0,\mu_0)^2\le 4\mu_0(|\nn L_0^{-1} (f-1)|^2),\ \ f\mu_0\in \scr P_0,$$
we obtain
$$ t\E^\nu\big[\W_2(\mu_{t,r,r}, \mu_0)^2\big|T<\tau\big] \le c \big\{r^{-\ff{d-4} 2} + 1_{\{d=4\}} \log r^{-1} + t^{-1} r^{-k}\big\},\ \ T\ge t\ge 1, r\in (0,1).$$
By this  and Lemma \ref{LL5}, we find a  decreasing function $c: (\ff{2(d-4)}3,\infty)\to (0,\infty)$ such that 
\beq\label{EW} \beg{split} &\E^\nu\big[\W_2(\mu_{t}, \mu_0)^2\big|T<\tau\big] 
\le c(k)  \big\{t^{-1}r^{-\ff{d-4}2} +t^{-1} 1_{\{d=4\}} \log r^{-1} + t^{-2} r^{-k}+r\big\},\\
&\ T\ge t\ge 1, r\in (0,1), k>\ff{2(d-4)}3.\end{split}\end{equation} 

(a) Let $d=4$. We take $r=t^{-1}$ for $t>1$, such that \eqref{EW} implies \eqref{R3}   for some constant $c>0$. 

(b) When $d\ge 5$.  Since 
$$\lim_{k\downarrow \ff{2(d-4)} 3} \Big\{2- \ff{2k}{d-2}\Big\} = \ff{2d+4}{3(d-2)}  > \ff 2 {d-2},$$
there exists $k> \ff{2(d-4)} 3$ such that $2- \ff{2k}{d-2}> \ff 2 {d-2}.$ 
So, we may take   $r= t^{-\ff 2 {d-2}}$ for $t>1$  such that \eqref{EW} implies the inequality   in  (3). 
\end{proof} 

\section{Lower bound estimate}
This section devotes to the proof of the following result, which together with Proposition \ref{Pn1} implies Theorem \ref{T1.1}. 

\beg{prp}\label{PN2} Let  $\nu\in \scr P_0$. There exists a constant $c>0$ such that \eqref{R00} holds, and when $\pp M$ is convex it holds for $c=1.$ 
   Moreover,  when $d\ge 5$, there exists a constant $c'>0$ such that 
 \beq\label{R5}\inf_{T\ge t} \big\{t\E[\W_2(\mu_t,\mu_0)|T<\tau]\big\}\ge c' t^{-\ff 2{d-2}},\ \ t\ge 1.\end{equation} 
\end{prp}

To estimate the Wasserstein distance from below, we use the idea of \cite{AMB} to construct a pair of functions in Kantorovich's dual formula,   which leads to the following lemma.  

\beg{lem}\label{LF2} There exists a constant $c>0$ such that
$$\W_2(\mu_{t,r},\mu_0)^2\ge \mu_0(|\nn L_0^{-1}(\rr_{t,r}-1)|^2) -c \|\rr_{t,r}-1\|_\infty^{\ff 7 3}(1+\|\rr_{t,r}-1\|_\infty^{\ff 1 3}),\ \ t,r>0.$$
\end{lem} 
\beg{proof} Let $f= L_0^{-1}(\rr_{t,r}-1)$, and take 
$$\varphi_\theta^\vv=-\vv\log P_{\ff{\vv \theta}2}^0 \e^{-\vv^{-1}f},\ \ \theta\in [0,1], \vv>0.$$
We have $\varphi_0=f$ and by \cite[Lemma 2.9]{WZ20}, 
\beg{align*}&\varphi_1^\vv(y)-f(x)\le \ff 1 2\big\{\rr(x,y)^2 +\vv \|(L_0f)^+\|_\infty +c_1 \vv^{\ff 1 2} \|\nn f\|_\infty^2\big\},\\
&\mu_0(f-\varphi_1^\vv)\le \ff 1 2 \mu_0(|\nn f|^2) + c_1 \vv^{-1} \|\nn f\|_\infty^4.\end{align*}
Since $L_0f=\rr_{t,r}-1$, this and the integration by parts formula imply
\beq\label{MF1} \beg{split}& \ff 1 2 \W_2(\mu_{t,r},\mu_0)^2 + \vv \|\rr_{t,r}-1\|_\infty + c_1\vv^{\ff 1 2 }\|\nn f\|_\infty^2 
\ge \mu_0(\varphi_1^\vv) - \mu_{t,r}(f)\\
& = \mu_0(\varphi_1^\vv-f) -\mu_0(f L_0 f)
\ge \ff 1 2 \mu_0(|\nn L_0^{-1} f|^2) - c_1\vv^{-1} \|\nn f\|_\infty^4,\ \ \vv>0.\end{split}\end{equation}
Next, by Lemma \ref{L*}(1) for $p=\infty$ and \eqref{000}, we find   constants $c_2,c_3,c_4>0$ such that 
\beg{align*} &\|\nn f\|_\infty =\|\nn L_0^{-1} (\rr_{t,r}-1)\|_\infty\le \int_0^\infty \|\nn P_s^0(\rr_{t,r}-1)\|_\infty\d s\\
&\le c_2 \int_0^\infty (1+s^{-\ff 1 2}) \|P_{s/2}^0 (\rr_{t,r}-1)\|_\infty\d s\\
&\le c_3 \|\rr_{t.r}-1\|_\infty \int_0^\infty (1+s^{-\ff 1 2})\e^{-(\ll_1-\ll_0)s/2} \d s \le c_4  \|\rr_{t.r}-1\|_\infty.\end{align*}
Combining this with \eqref{MF1} we find a constant $c_5>0$ such that 
$$ \W_2(\mu_{t,r},\mu_0)^2\ge \mu_0(|\nn L_0^{-1} f|^2)-c_5 \big\{\vv \|\rr_{t,r}-1\|_\infty +  \vv^{\ff 1 2} \|\rr_{t,r}-1\|_\infty^2 + \vv^{-1} \|\rr_{t,r}-1\|_\infty^4\big\},\ \ \vv>0.$$
By taking $\vv= \|\rr_{t,r}-1\|_\infty^{\ff 4 3}$ we finish the proof. 
\end{proof}
By Lemma \ref{LF2}, to derive a sharp lower bound of   $\W_2(\mu_{t,r},\mu_0)^2$,   we need to estimate $\|\rr_{t,r}-1\|_\infty$ and $\E^\nu\big [ \mu_0(|\nn L_0^{-1} (\rr_{t,r}-1)|^2)\big|T<\tau\big],$ which are included in the following three lemmas.  

\beg{lem}\label{LF1} For any $r>0$ and $\nu=h\mu$ with $\|h\phi_0^{-1}\|_\infty<\infty$, there exists a constant $c(r)>0$ such that 
$$\sup_{T\ge t} \E^{\nu}\big[\|\rr_{t,r}-1\|_\infty^{4} \big|T<\tau\big] \le c(r) t^{-2},\ \ t\ge 1.$$\end{lem}
\beg{proof} By Lemma \ref{LEM} and \eqref{EQ0}, it suffices to prove for $\nu=\nu_0$ and $T=t$ replacing $T\ge t$, i.e. 
for  a constant $c(r)>0$  we have 
\beq\label{LLM}   \E^{\nu_0}\big[\|\rr_{t,r}-1\|_\infty^{4} \big|t<\tau\big] \le c(r) t^{-2},\ \ t\ge 1.\end{equation}
By \eqref{L009}, \eqref{PR1}, \eqref{000}, and $\|\phi_0^{-1}\|_{L^2(\mu_0)}=1$, we find a constant $c_1>0$ such that 
\beg{align*}&\E^{\nu_0} [1_{\{r_1<\tau\}} |g_s(r_1,r_2)|^2]\\
& \le c_1 \e^{-\ll_0(2s-\ll_1)} \|\hat\phi_m\|_\infty^4 \big\{\e^{-(\ll_m-\ll_0)(r_2-r_1)}+\e^{-(\ll_1-\ll_0)(s-r_2)}\big\},\ \ s>r_2>r_1>0.\end{align*} 
By \eqref{QK0'} and $\P^{\nu_0}(t<\tau)\ge c_0 \e^{-\ll_0 t}$ for some constant $c_0>0$ and all $t\ge 1$,  this implies 
$$\E^{\nu_0} [ |\psi_m(t)|^4|t<\tau ] :=\ff{\E^{\nu_0} [ |\psi_m(t)|^41_{\{t<\tau\}} ] }{P^{\nu_0}(t<\tau)}  \le c_2 \|\hat\phi_m\|_\infty^4t^{-2},\ \ m\ge 1, t>1$$ for some constant $c_2>0$. Combining with  \eqref{DS1} gives  
\beg{align*} &   \E^{\nu_0}\big[\|\rr_{t,r}-1\|_\infty^{4} \big|t<\tau\big] \\
&\le  \bigg(\sum_{m=1}^\infty \e^{-(\ll_m-\ll_0)r} \|\hat\phi_m\|_\infty^{\ff 4 3}\bigg)^{3} \sum_{m=1}^\infty \e^{-(\ll_m-\ll_0)r} e^{\ll_0 t} \E^{\nu_0} [1_{\{r_1<\tau\}}|\psi_m(t)|^4]\\
&\le  \bigg(\sum_{m=1}^\infty \e^{-(\ll_m-\ll_0)r} \|\hat\phi_m\|_\infty^{\ff 4 3}\bigg)^{3} c_2 t^{-2} \sum_{m=1}^\infty \e^{-(\ll_m-\ll_0)r} \|\hat\phi_m\|_\infty^4. \end{align*}
By  \eqref{EG} and \eqref{CC},  this implies   \eqref{LLM}   for some constant $c(r)>0.$
\end{proof}

\beg{lem}\label{LF3}  Let $\nu=h\mu$ with $\|h\phi_0^{-1}\|_\infty<\infty$. Then for any $r>0$ there exists a constant $c(r)>0$ such that 
$$ \sup_{T\ge t}  \Big|t  \E^\nu\big [ \mu_0(|\nn L_0^{-1} (\rr_{t,r}-1)|^2)\big|T<\tau\big]- 2\sum_{m=1}^\infty \ff{\e^{-2(\ll_m-\ll_0) r } }{(\ll_m-\ll_0)^2}\Big|
\le \ff{c(r)} t, t\ge 1. $$\end{lem} 
\beg{proof}Let  $\{J_i:i=1,2,3\}$ be in \eqref{SD7}. By \eqref{000}, \eqref{CC}, and $ \|\hat\phi_m\|_{L^2(\mu_0)}=1$, we find a constant  $c_1>0$ such that   for any $T\ge t\ge s_2\ge s_1>0,$   
\beg{align*} |J_1(s_1,s_2)|& \le \|h\phi_0^{-1}\|_{\infty} \|P_{s_1}^0-\mu_0\|_{L^\infty(\mu_0)} \|\phi_m\phi_0^{-1}|_\infty^2 \|P_{T-s_2}^0-\mu_0\|_{L^1(\mu_0)} \|\phi_0^{-1}\|_{L^1(\mu_0)}\\
&  \le c_1\|\phi_m\phi_0^{-1}\|_\infty^2 \e^{-(\ll_1-\ll_0)(t+s_1-s_2)},\\ 
 |J_2(s_1,s_2)| &\le \|\phi_0\|_\infty \e^{-(\ll_m-\ll_0)(s_2-s_1)} \|h\phi_0^{-1}\|_\infty \|P_{s_1}^0-\mu_0\|_{L^\infty (\mu_0)} 
 \\
& \le c_1 \e^{-(\ll_1-\ll_0)s_2},\\
 |J_3(s_1,s_2)| &\le \|\phi_0\|_\infty \e^{-(\ll_m-\ll_0)(s_2-s_1)}  \|\phi_m\phi_0^{-1}\|_\infty^2\|P_{T-s_2}^0-\mu_0\|_{L^1(\mu_0)}\|\phi_0^{-1}\|_{L^1(\mu_0)}  \\
&  \le c_1 \|\phi_m\phi_0^{-1}\|_\infty^2  \e^{-(\ll_1-\ll_0)(t-s_1)}.\end{align*}
Substituting these into \eqref{007} and applying \eqref{EG} and \eqref{CC}, we find a constant $c(r)>0$ such that the desired estimate holds. 
\end{proof}

\beg{lem}\label{LF4} Let $\nu=h\mu$ with $\|h\phi_0^{-1}\|_\infty<\infty$. Then for any $r>0$  and $p\ge 2$, there exists a constant $c(r,p)>0$ such that 
$$\|\nn L_0^{-1} (\rr_{t,r}-1)|^{2p} \|_{L^{2p}(\mu_0)}  \le c(r,p),\ \ t>0.$$ 
\end{lem}
\beg{proof} Since $\rr_{t,r} =\ff 1 t\int_0^t p_r^0(X_s,\cdot) \d s,$ we have $\mu_0(\rr_{t,r})=1$ and $\|\rr_{t,r}\|_\infty\le \|p_r^0\|_\infty  <\infty$. 
Then by \eqref{000} and $\|\phi_0^{-1}\|_{L^2(\mu_0)}=1$,  we find  a constant $c_1(r)>0$ such that 
\beg{align*}&\mu_0\big(\phi_0^{-1}  \{P_{\ff s 4}^0|P_{\ff{3s}4}^0(\rr_{t,r}-1)|^p\}^2\big) \le \|\phi_0^{-1}\|_{L^2(\mu_0)}\|(P_{\ff{3s}4}^0-\mu_0)\rr_{t,r}\|_{L^{4p}(\mu_0)}^{2p}\\
&\le \|P_{\ff{3s}4}^0-\mu_0\|_{L^{4p}(\mu_0)}^{2p}\|\rr_{t,r}\|_\infty^{2p}\le c_1(r) \e^{-3(\ll_1-\ll_0)s}.\end{align*}
Combining this with \eqref{MFN} for $\vv=1$ and $\theta\in (0, \ff 1 {\ll_1-\ll_0})$, we finish the proof. 

\end{proof} 

Finally, since $\mu_{t,r}=\mu_t P_r^0$, to derive a lower bound of  $\W_2(\mu_t,\mu_0)$ from that of $\W_2(\mu_{t,r},\mu_0)$, we present the following result. 

\beg{lem}\label{LF5} There exist two constants $K_1,K_2>0$ such that for any probability measures $\mu_1,\mu_2 $ on $M^\circ$,
\beq\label{ND} \W_2(\mu_1P_t^0, \mu_2P_t^0)\le K_1\e^{K_2 t} \W_2(\mu_1,\mu_2),\ \ t\ge 0.\end{equation} 
When $\pp M$ is convex, this estimate holds for $K_1=1.$ 
\end{lem}
\beg{proof} When $\pp M$ is convex, by \cite[Lemma 2.16]{W20},   there exists a  constant $K$ such that 
$$\Ric- \Hess_{V+2\log\phi_0}\ge -K,$$   so that the desired estimate holds for
$K_1=1$ and $K_2=K$, see \cite{RS}.  

 In general, 
  following the line of \cite{W07}, we  make the boundary from non-convex to convex by using a conformal
change of metric.
Let  $N$ be the inward normal unit vector field of $\pp M$. Then the second fundamental form of $\pp M$
is a two-tensor on the tangent space of $\pp M$ defined by
$$\II(X,Y):=-\<\nn_X N, Y\>,\ \ X,Y\in T\pp M.$$  
Since $M$ is compact, we find a function $f\in C_b^\infty(M)$ such that $f\ge 1, N\parallel \nn f$ on $\pp M$,  
and $N\log f|_{\pp M} +\II(u,u) \ge 0$ holds on $\pp M$ for any $u\in T\pp M$ with $|u|=1.$ 
By \cite[Lemma 2.1]{W07} or \cite[Theorem 1.2.5]{W14}, $\pp M$ is convex under the metric
$$\<\cdot,\cdot\>'= f^{-2} \<\cdot,\cdot\>.$$ Let $\DD'$,   $\nn'$ and $\Hess'$ be the Laplacian, gradient  and Hessian induced
by the new metric $\<\cdot,\cdot\>'$. 
We have $\nn'= f^2 \nn$ and  (see (2.2) in \cite{TW})
$$L_0= f^{-2} \DD'+ f^{-2} \nn'\{V+2\log\phi_0+(d-2)f^{-1}\}.$$ Then the $L_0$-diffusion process $X_t^0$ with $X_0^0$ having distribution $\mu_1$ can be constructed by solving the 
following It\^o SDE on $M^\circ$ with metric $\<\cdot,\cdot\>'$ (see \cite{ATW06})
\beq\label{AE1} \d^I X_t^0= \big\{ f^{-2} \nn'(V+2\log\phi_0+(d-2)f^{-1})\big\}(X_t^0) \d t + \ss 2 f^{-1} (X_t^0)U_t\d B_t,\end{equation} 
where  $B_t$ is the $d$-dimensional Brownian motion, and $U_t$ is the horizontal lift of $X_t^0$ to the frame bundle $O'(M)$ with respect to  the metric $\<\cdot,\cdot\>'$.

Let $Y_0^0$ be a random variable independent of $B_t$ with distribution $\mu_2$ such that  
\beq\label{O1} \W_2(\mu_1,\mu_2)^2 = \E[\rr(X_0^0,Y_0^0)^2].\end{equation} 
For any $x,y\in M^\circ$, let $P'_{x,y}: T_xM\to T_y M$ be the parallel transform along the minimal geodesic from $x$ to $y$  induced by the metric $\<\cdot,\cdot\>'$, which is contained in $M^\circ$ by the convexity.   
Consider the coupling by parallel displacement  
\beq\label{AE2} 
\d^I Y_t^0= \big\{ f^{-2} \nn'(V+2\log\phi_0+(d-2)f^{-1})\big\}(Y_t^0) \d t + \ss 2 f^{-1} (Y_t^0) P_{X_t^0,Y_t^0}'U_t\d B_t.\end{equation}
  As explained in \cite[Section 3]{ATW06}, we may assume that 
$(M^\circ, \<\cdot,\cdot\>') $ does not have cut-locus such that $P'_{x,y}$ is a smooth map, which  ensures the existence and uniqueness of  $Y_t^0$. 
Since the distributions of $X_0^0$ and $Y_0^0$ are $\mu_1,\mu_2$ respectively, the law of   $(X_t^0,Y_t^0)$ is in the class $\C(\mu_1P_t^0,\mu_2P_t^0)$,  so that
\beq\label{JNC}   \W_2(\mu_1P_t^0, \mu_2P_t^0)^2\le \E[\rr(X_t^0,Y_t^0)^2],\ \ t\ge 0.\end{equation} 
Let $\rr'(x,y)$ be  the Riemannian distance between  $x$ and $y$ induced by $\<\cdot,\cdot\>':=f^{-2}\<\cdot,\cdot\>.$ By $1\le f\in C_b^\infty(M)$ we have 
\beq\label{MTT} \|f\|_\infty^{-1} \rr\le\rr'\le\rr.\end{equation}
Since except the term $f^{-2} \nn'\log \phi_0$, all coefficients in the SDEs are in $C_b^\infty(M)$, by It\^o's formula, there exists a constant $K$ such that 
\beq\label{JNN} \d \rr'(X_t^0, Y_t^0)^2 \le \big\{K\rr'(X_t^0,Y_t^0)^2 + I \big\}\d t + \d M_t,\end{equation}
where $M_t$ is a martingale and
\beg{align*} I:= \<(f^{-2} \nn'\log\phi_0)(\gg_1), \dot \gg_1\>'-  \<(f^{-2} \nn'\log\phi_0)(\gg_0), \dot \gg_0\>'.\end{align*}
Let $\gg: [0, 1]\to M$ be  the minimal geodesic  from $X_t^0$ to $Y_t^0$ induced by the metric $\<\cdot,\cdot\>'$,   which is contained in $M^\circ$ by the convexity, we obtain 
\beg{align*} I&= \int_0^1 \ff{\d}{\d s}  \<(f^{-2} \nn'\log\phi_0)(\gg_s), \dot \gg_s\>' \d s \\
&= \int_0^1\Big\{ \ff{f^{-2}(\gg_s)\Hess'_{\phi_0}(\dot\gg_s,\dot\gg_s)+\<\nn' f^{-2}(\gg_s),\dot \gg_s\>'\<\nn'\phi_0(\gg_s), \dot\gg_s\>'}{\phi_0(\gg_s)}
 - \ff{\{\<\nn'\phi_0(\gg_s),\dot\gg_s\>'\}^2}{(f^2\phi_0^2)(\gg_s)}\Big\}\d s\\
 &\le  \int_0^1\Big\{ (\phi_0^{-1}f^{-2})(\gg_s)\Hess'_{\phi_0}(\dot\gg_s,\dot\gg_s)+\ff {f^2} 4    \big[\<\nn' f^{-2}(\gg_s), \dot\gg_s\>'\big]^2\Big\}\d s\le C\rr'(X_t^0,Y_t^0)^2\end{align*} for some constant $C>0$,
 where the last step is due to $\<\dot\gg_s,\dot\gg_s\>'=\rr'(X_t^0,Y_t^0)^2$, $1\le f\in C_b^\infty(M)$, and that by the proof of \cite[Lemma 2.1]{W20} the convexity of $\pp M$ under $\<\cdot,\cdot\>'$ implies 
 $\Hess'_{\phi_0}\le c \phi_0$ for some constant $c>0$.  This and \eqref{JNN} yield
 $$\E[ \rr'(X_t^0, Y_t^0)^2]\le \E[ \rr'(X_0^0, Y_0^0)^2]\e^{(K+C)t},\ \ t\ge 0.$$
 Combining this with \eqref{O1} and \eqref{MTT}, we prove   \eqref{ND} for some constant $K_1,K_2>0.$
\end{proof} 

We are now ready to prove the main result in this section. 

\beg{proof}[Proof of Proposition \ref{PN2}]  (a) According to \eqref{NL1}, it suffices to prove 
 for $\nu=h\mu$ with $\|h\phi_0^{-1}\|_\infty<\infty$. 
Let $r>0$ be fixed. By Lemma \ref{LF2},  we obtain
\beq\label{JJ1} \beg{split} &t\E^\nu\big[\W_2(\mu_{t,r},\mu_0)^2\big|T<\tau\big] \ge t\E^\nu\big[1_{\{\|\rr_{t,r}-1\|_\infty \le\vv\}}\W_2(\mu_{t,r},\mu_0)^2\big|T<\tau\big] \\
&\ge t\E^\nu\big[1_{\{\|\rr_{t,r}-1\|_\infty \le\vv\}}\mu_0(|\nn L_0^{-1}(\rr_{t,r}-1)|^2) \big|T<\tau\big]-c\vv^2\\
&\ge t\E^\nu\big[ \mu_0(|\nn L_0^{-1}(\rr_{t,r}-1)|^2) \big|T<\tau\big]-c\vv^2\\
&\qquad -t\E^\nu\big[1_{\{\|\rr_{t,r}-1\|_\infty >\vv\}}\mu_0(|\nn L_0^{-1}(\rr_{t,r}-1)|^2) \big|T<\tau\big] ,\ \ \vv>0, T\ge t.\end{split}\end{equation}
By Lemma \ref{LF1} and Lemma \ref{LF4} with $p=3$, we find some constants $c_1,c_2>0$ such that 
\beg{align*} & t\E^\nu\big[1_{\{\|\rr_{t,r}-1\|_\infty >\vv\}}\mu_0(|\nn L_0^{-1}(\rr_{t,r}-1)|^2) \big|T<\tau\big] \le c_1t \big\{\P^\nu\big(\|\rr_{t,r}-1\|_\infty>\vv\big| T<\tau\big)\big\}^{\ff 2 3}\\
&\le c_1t \vv^{-\ff 8 3}\big\{ \E^\nu\big(\|\rr_{t,r}-1\|_\infty^4\big| T<\tau\big)\big\}^{\ff 2 3} \le c_2 \vv^{-\ff 8 3} t^{-\ff 1 3}, \ \ T\ge t.\end{align*}
Combining this with \eqref{JJ1} and Lemma \ref{LF3}, we find a constant $c_3>0$ such that 
\beg{align*}& t\E^\nu\big[\W_2(\mu_{t,r},\mu_0)^2\big|T<\tau\big] \ge t\E^\nu\big[\mu_0(|\nn L_0^{-1}(\rr_{t,r}-1)|^2)\big|T<\tau\big] -\vv_t\\
& \ge 2\sum_{m=1}^\infty \ff {\e^{-2(\ll_m-\ll_0)r}}{(\ll_m-\ll_0)^2} -\vv_t - c_3 t^{-1},\ \ T\ge t\ge 1,\end{align*} where
$$\vv_t:=\inf_{\vv>0} \{c\vv^2 + c_2 \vv^{-\ff 8 3} t^{-\ff 1 3}\}\to 0\ \text{as}\ t\to\infty.$$
Therefore, 
$$ \liminf_{t\to\infty} \inf_{T\ge t} \Big\{tE^\nu\big[\W_2(\mu_{t,r},\mu_0)^2\big|T<\tau\big] \Big\}\ge 2\sum_{m=1}^\infty \ff {\e^{-2(\ll_m-\ll_0)r}}{(\ll_m-\ll_0)^2},\ \ r>0.$$
Combining this with   Lemma \ref{LF5},  we derive 
$$ \liminf_{t\to\infty} \inf_{T\ge t} \Big\{tE^\nu\big[\W_2(\mu_{t},\mu_0)^2\big|T<\tau\big] \Big\}
\ge 2K_1^{-1} \e^{-K_1r} \sum_{m=1}^\infty \ff {\e^{-2(\ll_m-\ll_0)r}}{(\ll_m-\ll_0)^2},\ \ r>0.$$
Letting $r\to 0$ we prove \eqref{R00} for $c= K_1^{-1}$. By Lemma \ref{LF5}, we may take $c=1$ when $\pp M$ is convex. 

(b) The second assertion can be proved as in \cite[Subsection 4.2]{WZ20}.   For any $t\ge 1$ and $N\in\mathbb N$, let $\mu_{N}:= \ff 1 N\sum_{i=1}^N \dd_{X_{t_i}},$ where $t_i:=\ff{(i-1)t}N, 1\le i\le N.$ \cite[Proposition 4.2]{RE1} (see also \cite[Corollary 12.14]{RE2}) implies
\beq\label{JX0} \W_1(\mu_N,\mu_0)^2\ge c_0N^{-\ff 2 d},\ \ N\in\mathbb N, t\ge 1\end{equation} 
for some constant $c_0>0$. 
Write
$$\mu_t= \ff 1 N \sum_{i=1}^N \ff Nt \int_{t_i}^{t_{i+1}}\dd_{X_{s}}\d s.$$
By the convexity of $\W_2^2$, which follows from the Kantorovich dual formula, we have
\beq\label{JX00} \W_2(\mu_N,\mu_t)^2\le \ff 1 N \sum_{i=1}^N  \ff N t \int_{t_i}^{t_{i+1}} \W_2(\dd_{X_{t_i}}, \dd_{X_s})^2\d s 
= \ff 1 t \sum_{i=1}^N   \int_{t_i}^{t_{i+1}} \rr(X_{t_i},  X_s)^2\d s\end{equation} 

On the other hand,  by the Markov property, 
\beq\label{LX1} \E^\nu[\rr(X_{t_i}, X_s)^21_{\{T<\tau\}}] = \E^\nu\big[1_{\{t_i<\tau\}} P_{s-t_i}^D\{\rr(X_{t_i},\cdot)^2 P_{T-s}^D 1\}(X_{t_i})\big].\end{equation}
Since $P_t^D 1\le c_1\e^{-\ll_0t}$ for some constant $c_1>0$ and all $t\ge 0$,  \eqref{PR0} implies
\beq\label{LX2}\beg{split} & P_{s-t_i}^D \{\rr(x,\cdot)^2 P_{T-s}^D 1\}(x) \\
&\le c_1 \e^{-\ll_0(T-s)} P_{s-t_i}^D \rr(x,\cdot)^2(x) \le c_1  \e^{-\ll_0(T-s)}  \phi_0 (x) P_{s-t_i}^0\{\rr(x,\cdot)^2\phi_0^{-1}\}(x).\end{split} \end{equation}
It is easy to see that 
$$L_0\{\rr(x,\cdot)^2\phi_0^{-1}\}\le c_2\phi_0^{-2}$$ holds on $M^\circ$ for some constant $c_2>0$. So, by \eqref{*B*},  we find  a constant $c_3>0$ such that 
$$P_{s-t_i}^0\{\rr(x,\cdot)^2\phi_0^{-1}\}(x)\le c_2 \E^x\int_0^{s-t_i} \phi_0^{-2}(X_r)\d r\le c_3 (s-t_i) \log (1+\phi_0^{-1}(x)).$$
Combining this with \eqref{LX1} and \eqref{LX2}, and using $P_t^D 1\le c_1\e^{-\ll_0t}$ observed above, we find a constant $c_5>0$ such that 
\beg{align*} & \E^\nu[\rr(X_{t_i}, X_s)^21_{\{T<\tau\}}] \le c_4 \e^{-\ll_0 T} \nu(\log (1+\phi_0^{-1})) (s-t_i) \\
&\le c_4\|h\phi_0^{-1}\|_\infty \mu(\phi_0 \log (1+\log \phi_0^{-1}))(s-t_i)\e^{-\ll_0 T} \le 
c_5 (s-t_i)\e^{-\ll_0 T},\ \ s\ge t_i.\end{align*} 
Since $\P^\nu(T<\tau)\ge c_0\e^{-\ll_0 T}$ for some constant $c_0>0$ and all $T\ge 1$, we find a constant $c>0$ such that 
$$\E^\nu[\rr(X_{t_i}, X_s)^2|T<\tau] \le c(s-t_i),\ \ s\ge t_i.$$
Combining this with \eqref{JX0} and \eqref{JX00}, we find a constant $c_6>0$ such that 
$$\E^\nu[\W_1(\mu_t,\mu_0)^2|T<\tau]\ge \ff {c_1}  2 N^{-\ff 2 d}- c_6t N^{-1},\ \ T\ge t.$$
Taking $N=\sup\{i\in\mathbb N: i\le \aa t^{\ff d{d-2}}\}$ for some $\aa>0$, we derive
$$t^{\ff 2 {d-2}}     \inf_{T\ge t}  \{\E^\nu [\W_1(\mu_0,\mu_t)^2|T<\tau] \} \ge  \ff {c_2}{2\aa^{\ff 2 d}}- \ff {2c'}\aa,\ \ t\ge 1.$$   Therefore,
$$   t^{\ff 2 {d-2}}  \inf_{T\ge t}   \E^\nu [\W_1(\mu_0,\mu_t)^2 |T<\tau]\ge \sup_{\aa>0} \Big(\ff {c_2}{2\aa^{\ff 2 d}}- \ff {2c'}\aa\Big)>0,\ \ t\ge 1.$$

\end{proof}

\end{document}